\documentclass{article}[11pt,titlepage]

\newtheorem{theorem}{Theorem}[section]

\newcommand\nc{\newcommand}
\nc{\od}[2]{\frac{d#1}{d#2}}
\nc{\be}{\begin{equation}}
\nc{\ee}{\end{equation}}
\nc{\bd}{\begin{displaymath}}
\nc{\ed}{\end{displaymath}}
\nc{\bea}{\begin{eqnarray}}
\nc{\eea}{\end{eqnarray}}
\nc{\p}{\partial}

\usepackage{graphics}
\usepackage{amsmath}
\usepackage{amsfonts}
\usepackage{amssymb}

\begin{document}
\title{Multi-phase Stefan problems for a nonlinear 1-d model of cell-to-cell adhesion and diffusion}
\author{K. Anguige\thanks{e-mail: keith.anguige@oeaw.ac.at} \\ {\em RICAM, Austrian Academy of Sciences,}\\{\em Altenbergerstr. 69, A-4040 Linz, Austria}}

\date{September 22nd, 2009}
\maketitle

\begin{abstract}
We consider a family of multi-phase Stefan problems for a certain 1-d model of cell-to-cell adhesion and diffusion, which takes the form of a nonlinear forward-backward parabolic equation. In each material phase the cell density stays either high or low, and phases are connected by jumps across an `unstable' interval. We develop an existence theory for such problems which allows for the annihilation of phases and the subsequent continuation of solutions. Stability results for the long-time behaviour of solutions are also obtained, and, where necessary, the analysis is complemented by numerical simulations.
\end{abstract}

\section{Introduction}

In this paper, we give further consideration to the 1-d continuum model for adhesion/diffusion of biological cells developed by Anguige and Schmeiser in \cite{anguige}, which took the form of the nonlinear diffusion equation
\be
\frac{\partial\rho}{\partial t} = \frac{\partial}{\partial x}\left(D(\rho)\frac{\partial\rho}{\partial x}\right),
\label{cont_rho}\ee
with quadratic diffusivity
\be
D(\rho)= 3\alpha\left(\rho-\frac{2}{3}\right)^2+1-\frac{4}{3}\alpha,
\label{D}\ee
for the scaled cell density $\rho(x,t)\in[0,1]$, and the adhesion coefficient $\alpha\in[0,1]$, the boundary condition being just $\frac{\p\rho}{\p x}=0$ at $x=0, 1$.

These equations were obtained as the formal continuum limit of the fundamental biased-random-walk model

\be
\frac{\partial \rho_i}{\partial t}=\mathcal{T}^+_{i-1}\rho_{i-1}+\mathcal{T}^-_{i+1}\rho_{i+1}-(\mathcal{T}^+_i+\mathcal{T}^-_i)\rho_i~,\label{walk}
\ee
with transitional probabilities
\be
\mathcal{T}^{\pm}_i = (1-\rho_{i\pm 1})(1-\alpha\rho_{i\mp 1})/h^2 ,\label{Ti} 
\ee
on a lattice of points $x_i=ih$, by taking Taylor expansions about $x_i$, and letting $h\rightarrow 0$. In the derivation of this equation, $h$ was interpreted as a (microscopic) measure of cell size.

We recall from \cite{anguige} that (\ref{cont_rho})-(\ref{D}) is globally well posed if $\alpha<\frac{3}{4}$. If, on the other hand, $\alpha>\frac{3}{4}$ then (\ref{cont_rho})-(\ref{D}) is ill-posed iff the initial density profile protrudes into the `unstable' interval
\be
I_{\alpha}= (\rho^{\flat}(\alpha), \rho^{\sharp}(\alpha)):= \left(\frac{2\alpha-\sqrt{\alpha(4\alpha-3)}}{3\alpha},\frac{2\alpha+\sqrt{\alpha(4\alpha-3)}}{3\alpha}\right)\subset[1/3,1],\label{rhointerval}
\ee
since in that case $D(\rho)$ is positive iff $\rho\notin I_{\alpha}$, and positivity is preserved by the Maximum Principle.

For completeness, note that in the borderline case $\alpha=\frac{3}{4}$, equation (\ref{cont_rho}) is just the porous-medium equation with quadratic diffusivity and possible change of sign about $\rho=\frac{2}{3}$. For initial data which stays away from $\rho=\frac{2}{3}$ (either above or below), (\ref{cont_rho}) is uniformly parabolic, and global existence of a smooth solution follows as for $\alpha<\frac{3}{4}$, while for degenerate initial data one is merely guaranteed a (unique) globally existing weak solution \cite{vazquez}.

The ill-posedness of (\ref{cont_rho}) for $\alpha>\frac{3}{4}$ is related to the presence of fine (wavelength $O(h)$) spatial oscillations, as well as plateau formation, in solutions of the discrete system (\ref{walk}), and the absence of a straightforward existence theory for (\ref{cont_rho}) leads one to ask just what model should be taken as a reasonable continuum limit of (\ref{walk}) in the high-adhesion regime.

One approach, and the one we shall adopt in this paper, is to circumvent the problem of ill-posedness by simply declaring that $\rho$-values in $I_{\alpha}$ are forbidden, and considering solutions to (\ref{cont_rho})-(\ref{D}) which may jump across $I_{\alpha}$ (possibly multiple times), but which are otherwise smooth. Mathematically, one is then dealing with a kind of (multi-phase) Stefan problem for the density $\rho(x,t)$ and the jump locations $s_i(t)$, such that the $s_i$ are dynamically determined by local conservation of mass, or, in other words, by the Rankine-Hugoniot condition. 

In \cite{anguige}, the analysis of (\ref{walk}) was aided by considering higher-order modifications of the leading-order equation (\ref{cont_rho}). One such $O(h^2)$-modification takes the form of the fourth-order PDE
\be \label{mod-equ}
\frac{\partial\rho}{\partial t} = \frac{\partial^2}{\partial x^2}\left( K(\rho) + h^2\left(\alpha\rho(\rho-1) \frac{\partial^2\rho}{\partial x^2} 
- \alpha\rho \left( \frac{\partial\rho}{\partial x}\right)^2 + \frac{1}{12} \frac{\p\rho}{\p t} \right)\right),
\ee
where the cubic $K(\rho)$ is a primitive for $D(\rho)$. This equation is rather similar to the viscous Cahn-Hilliard equation \cite{vii}, and is a regularisation of (\ref{cont_rho}) in the sense that it is (at least locally) well-posed on $S^1$, for each fixed value of the microscopic parameter $h$, and for all $\alpha<1$. Presumably, solutions continue to exist globally, as for Cahn-Hilliard, but a proof is currently lacking.

The steady-state equation for (\ref{mod-equ}) can, after a change of variables, be written as a Hamiltonian dynamical system, and amongst the solutions there is, for each $\alpha$, a unique heteroclinic cycle. These heteroclinic cycles correspond to (two-level) plateau solutions of (\ref{mod-equ}), are close (for small $h$) to square-wave weak solutions of (\ref{cont_rho}), and their critical points, denoted by $\rho_1(\alpha)$ and $\rho_2(\alpha)$, such that $\rho_1<\rho^{\flat}<\rho^{\sharp}<\rho_2$ and $K(\rho_1)=K(\rho_2)$, match very well the numerically-observed long-time plateau values in solutions of (\ref{walk}) (see \cite{anguige}). For these reasons, we demand in our Stefan-problem framework that any jumps across $I_\alpha$ should connect $\rho_1(\alpha)$ to $\rho_2(\alpha)$. Furthermore, in order to avoid the degeneracies at $\partial I_{\alpha}$, we require that the initial data satisfy $\rho<\rho^{\flat}$ in low-density phases, and $\rho>\rho^{\sharp}$ in high-density ones.

We emphasise that $\rho_1$ and $\rho_2$ are determined by the particular choice of microscopic model (\ref{walk})-(\ref{Ti}); other model choices are possible, and these will result in different $\rho$-values.

The paper is organised as follows. In Section 2, we develop a partial existence theory and perform a steady-state analysis for the simplest Stefan problem, namely, that for which there is only a single discontinuity in the density; solution behaviour is further clarified with the aid of several numerical simulations. In Section 3, we extend the analysis to the general multi-phase case, which, in particular, allows for the annihilation of phases via coalescence events. Finally, in the Appendices, we collect a number of results from classical parabolic theory which are used throughout the paper.

\section{The 1-jump problem}

We begin the analysis by considering the simplest possible case, in which there are just two phases, connected by a single jump from $\rho_1$ to $\rho_2$; this situation is depicted in Figure 1. For definiteness, we will assume that the low-density phase lies to the left, and the high-density phase to the right, of the discontinuity; the converse arrangement can, of course, also be treated.

\setcounter{totalnumber}{1}

\begin{figure}

\centering

\resizebox{5in}{5in}{\includegraphics{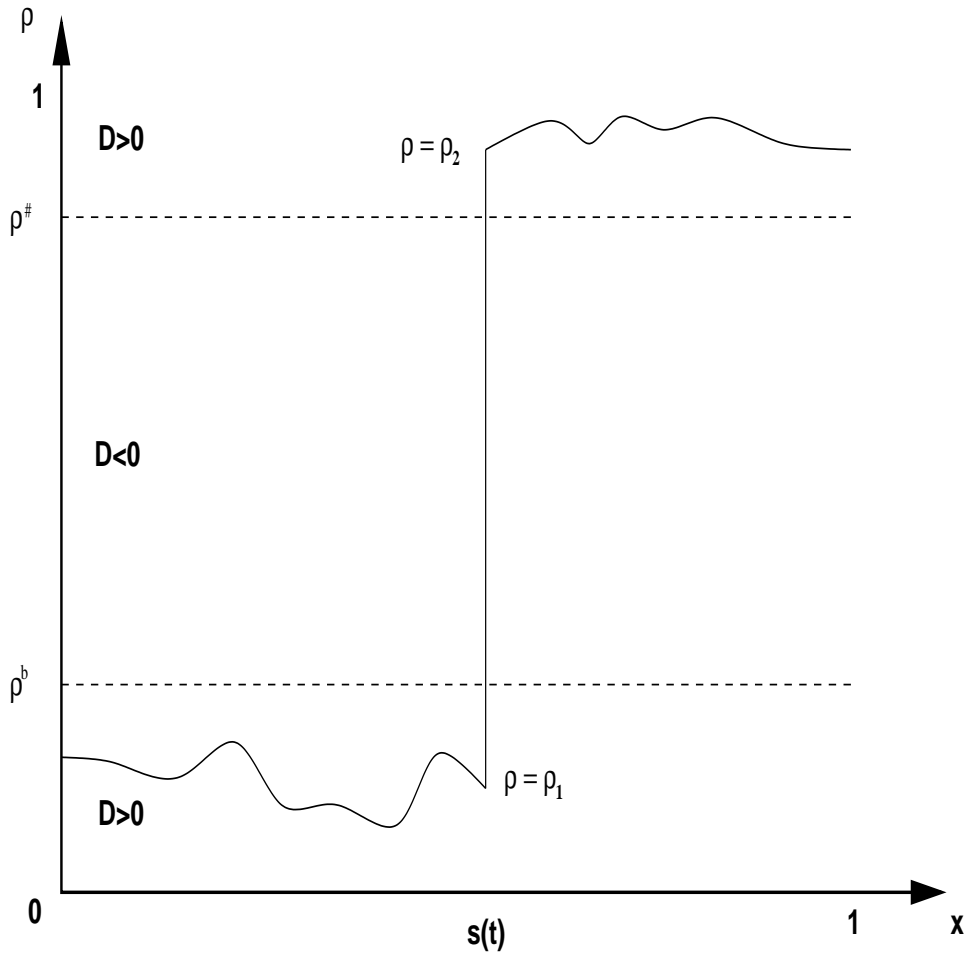}}

\caption{A typical density profile for the 1-jump Stefan problem, $SP_1$. }

\end{figure}

\subsection{Formulation}

Let the location of the jump discontinuity be denoted by $s(t)$. The Stefan problem for a given $\alpha>\frac{3}{4}$ then consists of looking for a function $\rho(x,t)$ on $[0,1]\times[0,T]$ which satisfies
\be
\frac{\partial\rho}{\partial t} = \frac{\partial}{\partial x}\left(D(\rho)\frac{\partial\rho}{\partial x}\right) ;\quad 0<x<s(t),\label{Lphase}g
\ee
subject to
\be
\frac{\p\rho}{\p x}(0,t)  =  0\quad,\quad\rho(s^-(t),t) =  \rho_1(\alpha),
\ee
and
\be
\frac{\partial\rho}{\partial t} =  \frac{\partial}{\partial x}\left(D(\rho)\frac{\partial\rho}{\partial x}\right) ;\quad s(t)<x<1,
\ee
subject to 
\be
\frac{\p\rho}{\p x}(1,t)  =  0\quad,\quad\rho(s^+(t),t) =  \rho_2(\alpha).\label{Rphase}
\ee

The evolution of $s(t)$ is determined by the Rankine-Hugoniot jump condition

\be
\frac{ds}{dt} = \frac{(J^+ - J^-)}{(\rho_2-\rho_1)},\label{RH}
\ee
where the fluxes are given by $J^{\pm}=-D(\rho(s^{\pm}(t),t))\p_x\rho(s^{\pm}(t),t)$, this equation being obtained by differentiating the statement of conservation of mass

\be
\int_0^{s(t)} \rho(x,t)~dx + \int_{s(t)}^1 \rho(x,t)~dx = M,\label{Mass}
\ee
and using (\ref{Lphase})-(\ref{Rphase}).

The initial datum for $\rho$, satisfying the boundary and jump conditions, is chosen to be smooth away from the initial discontinuity, with $\rho<\rho^{\flat}$ in the low-density phase and $\rho>\rho^{\sharp}$ in the high-density phase. By construction, $D(\rho)$ is then initially positive on each phase, and, by virtue of the Maximum Principle, we can reasonably expect the Stefan problem to be well posed. In the sequel, we will refer to the coupled system (\ref{Lphase})-(\ref{RH}) as $SP_1$; analogous problems $SP_n$, with $n$ jumps, will be treated in Section 3.

Finally, note that we allow for the possibility that $s(t)$ may hit the domain boundary at $x=0$ or $1$ in finite time. In this case, provided the gradient remains bounded as the boundary is approached, the solution can be continued via the ordinary Neumann problem for (\ref{cont_rho}), which we will call $NP$. As a convention, a solution continued in this way will still be referred to globally as a solution of $SP_1$.

\subsection{Steady states and their stability}

If the total mass $M$ satisfies $\rho_1<M<\rho_2$ then there is precisely one discontinuous steady-state solution of $SP_1$, given by the step function
\bea
\rho^{\ast}(x) & = &
\left\{
\begin{array}{ccc}
\rho_1 & : & x\in[0,s^{\ast})\\
\rho_2 & : & x\in(s^{\ast},1]
\end{array}
\right.
\nonumber\\
s^{\ast} & = &\frac{(\rho_2-M)}{(\rho_2-\rho_1)},\label{s_star}
\eea
while if $M\leq\rho_1$ or $M\geq\rho_2$ then no such discontinuous solution exists.

Moreover, $SP_1$ clearly has the uniform steady-state solution 
\be
\rho(x)=M,~0\leq x\leq1
\ee
if and only if $M\leq\rho^{\flat}$ or $M\geq\rho^{\sharp}$.

For the two possible kinds of steady state we have the following stability results:

\begin{theorem}
Suppose we have a global smooth solution pair, ($\rho(x,t), s(t)$), for $SP_1$, with initial data $(\rho_0(x),s_0)$, such that $0<s_0<1$. Then
\\~\\
(i) if $\rho^{\flat}<M<\rho^{\sharp}$, $\rho(x,t)$ converges exponentially to $\rho^{\ast}$ in $L^2$-norm, and $s(t)$ converges exponentially to $s^{\ast}$, as $t\rightarrow\infty$,
\\~\\
(ii) if $\rho_1<M\leq\rho^{\flat}$ (resp. $\rho^{\sharp}\leq M<\rho_2$), and $\rho_0\leq M~(\rho_0\geq M)$ in the low (high) phase, the solution is attracted towards $(\rho^{\ast},s^{\ast})$ as in (i),
\\~\\
(iii) if $M<\rho_1$ (resp. $M>\rho_2$), then $s(t)$ hits $x=0~(x=1)$ in finite time, and the continued solution converges to the uniform steady state $\rho=M$, exponentially in $L^2$.
\\~\\
(iv) if $M=\rho_1$ (resp. $M=\rho_2$), then $s\rightarrow 1$ $(0)$ as $t\nearrow T$, some $T$, with the possibility that $T=\infty$, and $\rho(x,t)$ approaches $\rho_1$ ($\rho_2$) exponentially in $L^2$ in the low (high) phase as $t\nearrow T$. If $T$ is finite, the continued solution of NP converges to the uniform steady state $\rho=M$, exponentially in $L^2$.
\end{theorem}

{\em Proof.}~~~~To prove (i), first note that, by (\ref{Mass}) and the fact that $\rho$ cannot enter $I_{\alpha}=(\rho^{\flat},\rho^{\sharp})$,
\be
0<s_{min}\leq s(t)\leq s_{max}<1~,\forall t,
\ee
where
\bea
s_{min}(M) & = & (\rho^{\sharp}-M)/\rho^{\sharp},\\
s_{max}(M) & = & (1-M)/(1-\rho^{\flat}).
\eea

Next note that, on the subinterval $(0,s(t))$, we have

\be
\frac{\p}{\p t}(\rho-\rho_1) = \frac{\p}{\p x}\left(D(\rho)\frac{\p}{\p x}(\rho-\rho_1)\right).
\ee
Hence, multiplying through by $(\rho-\rho_1)$, integrating by parts, and using the mixed boundary conditions gives

\be
\frac{1}{2}\frac{d}{dt}\int_0^{s(t)}(\rho-\rho_1)^2~dx  =  -\int_0^{s(t)}D(\rho)(\p_x(\rho-\rho_1))^2~dx,
\ee
which implies
\be
\frac{1}{2}\frac{d}{dt}\|\rho-\rho_1\|^2_{L^2(0,s)}  \leq  -\epsilon\|\p_x(\rho-\rho_1)\|^2_{L^2(0,s)},\label{L2}
\ee
where $\epsilon = \inf_{0<x<s(0)}\{D(\rho(x,0))\}>0$.

Since $\rho=\rho_1$ at $x=s^-(t)$, we also have the Poincar\'{e} inequality

\be\|\rho-\rho_1\|_{L^2(0,s)}\leq s\|\p_x(\rho-\rho_1)\|_{L^2(0,s)},
\ee
and hence (\ref{L2}) implies

\be
\|\rho-\rho_1\|^2_{L^2(0,s)}(t)\leq\|\rho-\rho_1\|^2_{L^2(0,s)}(0)e^{-\epsilon t/s^2_{max}}.
\ee
An analogous inequality on the subinterval $(s(t),1)$ is obtained in exactly the same way.

Next, from (\ref{Mass}) we have

\be
\int_0^{s(t)} (\rho(x,t)-\rho_1)~dx~+~\rho_1s(t)~+~\int_{s(t)}^1 (\rho(x,t)-\rho_2)~dx~+~\rho_2(1-s(t)) = M,
\ee
and hence, rearranging and using the $L^2$-decay just shown,

\be
|M-\rho_2-s(t)(\rho_1-\rho_2)|\leq c_1e^{-c_2t},\label{s_est1}
\ee
for some $c_1, c_2$.

Substituting the definition of $s^{\ast}$ from (\ref{s_star}) into (\ref{s_est1}) gives us

\be
|s-s^{\ast}|(t)\leq\frac{c_1e^{-c_2t}}{(\rho_2-\rho_1)},\label{s_est2}
\ee
as required.

To prove (ii), note, for example, that if $\rho_1<M\leq\rho^{\flat}$ and $\rho_0(x)\leq M$ for $0\leq x\leq s_0$ then $\rho=M$ is a supersolution in the low phase. From this, it follows once again that $s(t)$ remains bounded away from $0$ and $1$, by conservation of mass. Convergence to $(\rho^{\ast},s^{\ast})$ is proved as before. The same argument goes through for $\rho^{\sharp}\leq M<\rho_2$ and $\rho_0(x)\geq M$, $s_0\leq x\leq 1$.

For (iii), if $M<\rho_1$~(resp. $M>\rho_2$) then the only candidate for a steady-state solution is given by $\rho=M$ and $s=1~(s=0)$, and, by $L^2$-decay in the low (high) phase and conservation of mass, $s(t)$ is forced to hit the boundary in finite time; $L^2$-convergence for the subsequent $NP$ is proved via the usual energy estimate

Finally, if $M=\rho_1$ (resp. $\rho_2$), then the only possible steady state is $\rho=\rho_1$ $(\rho_2)$ for $0\leq x\leq1$. Clearly, $s(t)$ remains bounded away from $0$ ($1$) and therefore $\rho(\cdot,t)$ decays to $\rho_1$ ($\rho_2$), exponentially in $L^2$ on $(0,s)$ ($(s,1)$), as we have already seen. The $L^2$-decay implies that $s\rightarrow 1$ ($0$) as $t$ increases, but the convergence may take infinitely long in this exceptional case; if the convergence takes place in finite time then the subsequent $NP$ is as in (iii). This proves part (iv) $\square$
\\~\\
We conclude this discussion by noting that, for the cases not covered by Theorem 2.1, namely
\begin{enumerate}
 \item  $\rho_1<M\leq\rho^{\flat},~\rho_0(x)>M$~for some $x<s_0$, 
 \item  $\rho^{\sharp}\leq M<\rho_2,~\rho_0(x)>M$~for some $x>s_0$,
\end{enumerate}
we have not been able to find a clean analytical criterion for determining which of the two possible steady states will be approached at large times, for given initial data. However, in the next subsection we show that the discontinuous steady state (\ref{s_star}) is always linearly stable in a certain sense, whenever it exists, and, moreover, we report on numerical simulations which suggest that $s(t)$ can hit the domain boundary in finite time, provided condition 1 or 2 (directly above) is satisifed, and the initial density profile is sufficiently far from (\ref{s_star}); in other words, $SP_1$ appears to be bistable for some values of the total mass.

\subsection{Existence and uniqueness for smooth data}

\subsubsection{Background}
We can make a connection between our moving-boundary problem $SP_1$ and the traditional two-phase Stefan problem for the melting of ice in water by identifying $\rho$ as the specific latent heat, and $\sigma:=K(\rho)$ as the temperature. In the traditional formulation, our low-density phase is thought of as the solid, and our high-density phase as the liquid phase. Given this, one imposes the physically reasonable condition
\be
\rho_0\leq\rho_1~\textrm{in the solid (or low-density) phase,}\qquad \rho_0\geq\rho_2~\textrm{in the liquid (or high-density) phase},\label{class_data}
\ee
and then one can write down a favourable weak formulation of the problem, which can be shown to have a unique solution for bounded data (see Appendix C). Furthermore, this weak solution turns out to be a global-in-time smooth solution of the original problem (see, e.g., \cite{meir}).

For $SP_1$, however, the restriction on the initial data is just
\be
\rho_0<\rho^{\flat}~\textrm{in the low-density phase,}\qquad \rho_0>\rho^{\sharp}~\textrm{in the high-density phase},\label{SP_1_data}
\ee
which is weaker than (\ref{class_data}). In particular, this entails that $K(\rho)$ can no longer be assumed monotonically increasing, and, consequently, that there is no nice weak formulation of the problem. We note that, in the literature, uniqueness of solutions, given (\ref{class_data}), is usually proved using the weak formulation \cite{lady,meir}. Moreover, in the basic proof of local existence of smooth solutions as presented in \cite{meir}, for example, it is not so easy to discern whether the restriction (\ref{class_data}) can be relaxed to (\ref{SP_1_data}). For these reasons, we now present a concise, reasonably self-contained, classical existence-and-uniqueness theory for $SP_1$, with data subject merely to (\ref{SP_1_data}), which employs the standard H\"{o}lder and Sobolev estimates of linear parabolic theory.

\subsubsection{A local-existence theorem}

In order to solve the moving-boundary problem $SP_1$, we employ, in each phase, the standard trick of rescaling the spatial variable (see, for example, \cite{x}), such that (\ref{Lphase}) and (\ref{Rphase}) become a pair of fixed-boundary problems, coupled by the (rescaled) Rankine-Hugoniot condition.

Specifically, for $0\leq x\leq s(t)$ we introduce the coordinate transformation

\be
(x,t)\rightarrow(\hat{x},\hat{t}) :\quad\hat{x}=\frac{x}{s},\quad\hat{t}=t,
\ee
while for $s(t)\leq x\leq 1$ we take

\be
(x,t)\rightarrow(\tilde{x},\tilde{t}) : \quad\tilde{x}=\frac{x-s}{1-s},\quad\tilde{t}=t.
\ee

In terms of these new coordinates, and writing $v(\hat{x},\hat{t})=\rho(x,t)$ for $0\leq x\leq s(t)$, $w(\tilde{x},\tilde{t})=\rho(x,t)$ for $s(t)\leq x\leq 1$, equations (\ref{Lphase}) and (\ref{Rphase}) become, upon dropping hats,

\be
\frac{\p v}{\p t} = \frac{1}{s^2}\frac{\p}{\p x}\left(D(v)\frac{\p v}{\p x}\right) + x\frac{\dot{s}}{s}\frac{\p v}{\p x},\label{v}
\ee
for $x\in (0,1)$, subject to $v_x(0,t)=0,~v(1,t)=\rho_1$, and, respectively,

\be
\frac{\p w}{\p t} = \frac{1}{(1-s)^2}\frac{\p}{\p x}\left(D(w)\frac{\p w}{\p x}\right) + \dot{s}\frac{(1-x)}{(1-s)}\frac{\p w}{\p x},\label{w}
\ee
for $x\in (0,1)$, subject to $w(0,t)=\rho_2,~w_x(1,t)=0$, while the jump condition (\ref{RH}) assumes the form

\be
\frac{ds}{dt} = -\left(\frac{D(\rho_2)w_x(0,t)}{(1-s)} - \frac{D(\rho_1)v_x(1,t)}{s}\right)(\rho_2-\rho_1)^{-1},\label{rescaled_RH}
\ee
and we also have the rescaled conservation-of-mass equation 

\be
s(t) = \frac{M - \int_0^1 w~dx}{\left[\int_0^1 v~dx - \int_0^1 w~dx\right]}.\label{resc_COM}
\ee

Finally, for the application of classical parabolic theory, the initial data and the Dirichlet condition at $s_0:=s(0)$ must satisfy a certain first-order compatibility condition, which is obtained by setting $v_t$ and $w_t$ equal to zero in (\ref{v}) and (\ref{w}), and substituting the initial data, evaluated at the phase boundary, into the right-hand sides of (\ref{v}), (\ref{w}), and (\ref{rescaled_RH}).

The most compact way of writing this compatibility condition is to go back to the original coordinates and introduce the dependent (temperature) variable $\sigma=K(\rho)$, in terms of which the relevant equations are simply
\bea
-\sigma_x^+(\sigma_x^+ - \sigma_x^-) & = & D(\rho_1)\sigma_{xx}^+,\nonumber\\
-\sigma_x^-(\sigma_x^+ - \sigma_x^-) & = & D(\rho_2)\sigma_{xx}^-,
\eea
where, for example, $\sigma_x^{\pm}= \sigma_x(s^{\pm}(0),0)$.

We are able to prove the following:

\begin{theorem}
Given initial data $(v_0, w_0, s_0)$, such that $v_0, w_0\in C^{l+2}((0,1))$,~$l>0$,~$0<s_0<1$, $0\leq v_0<\rho^{\flat}$, $\rho^{\sharp}<w_0\leq 1$, and such that the first-order compatibility condition is satisfied, the system of equations (\ref{v}), (\ref{w}), (\ref{rescaled_RH}) (and hence problem $SP_1$) has a unique classical solution on some small time interval $[0,T]$.
\end{theorem}

{\em Proof.}~~~~The proof proceeds via an iterative scheme of successive approximations, and to get the required strong-convergence properties, we work in the setting of `parabolic' H\"{o}lder spaces.

First of all, we set $Q_T=[0,1]\times[0,T]$, and, as in \cite{lady}, \cite{meir}, for example, let $H^{l,l/2}({Q}_T)$, $l>0$ non-integer, denote the Banach space of functions $u$ with continuous derivatives $D^r_tD^s_xu$, for $2r+s\leq l$, equipped with the norm

\be
|u|^{(l)}_{Q_T} = \sum_{2r+s=0}^{[l]} \|D^r_tD^s_xu\|_{L^{\infty}(Q_T)} + \sum_{2r+s=[l]} \langle D^r_tD^s_xu\rangle^{(l-[l])}_{x,Q_T} + \sum_{0<l-2r-s<2}\langle D_t^rD_x^su\rangle^{(l-2r-s)/2}_{t,Q_T},
\ee
where, for $0<\beta<1$,

\bea
\langle v\rangle^{(\beta)}_{x,Q_T} & = & \sup_{(x,t),(x',t)\in Q_T}\left\{|v(x,t)-v(x',t)||x-x'|^{-\beta}\right\},\\
\langle v\rangle^{(\beta)}_{t,Q_T} & = & \sup_{(x,t),(x,t')\in Q_T}\left\{|v(x,t)-v(x,t')||t-t'|^{-\beta}\right\}.
\eea

Now we describe our iterative scheme for obtaining new approximate solutions of (\ref{v}), (\ref{w}), (\ref{rescaled_RH}) from old ones.

Suppose we have smooth $i$th iterates $(v^i,w^i)$, with positive diffusivities on $Q_T$, which satisfy the initial and boundary conditions, and suppose that 

\be
\max\left\{|v^i|^{(l+1)}_{Q_T}, |w^i|^{(l+1)}_{Q_T}\right\}\leq C,\label{vw_c_bound}
\ee
for some $l\in(0,1)$.

Next, determine the approximant $s^i(t)$ by solving the nonlinear ODE

\be
\frac{ds^i}{dt} = -\left(\frac{D(\rho_2)w^i_x(0,t)}{(1-s^i)} - \frac{D(\rho_1)v^i_x(1,t)}{s^i}\right)(\rho_2-\rho_1)^{-1},\label{iterated_RH}
\ee
which, for the given datum $s_0$, has a unique smooth solution on $[0,T]$, provided $T=T(C)$ is chosen small enough (Picard's Theorem). Clearly, the ordinary H\"{o}lder norm of $s^i$ is bounded according to 

\be
\|\dot{s}^i\|_{C^{l/2}([0,T])}\leq F(C),\label{sdot_bound} 
\ee
for some function $F$.

The functions $v^i, w^i, s^i$ are now inserted into the right-hand sides of the linearised field equations

\be
\frac{\p v^{i+1}}{\p t} = \frac{1}{(s^i)^2}\frac{\p}{\p x}\left(D(v^i)\frac{\p v^{i+1}}{\p x}\right) + x\frac{\dot{s}^i}{s^i}\frac{\p v^{i+1}}{\p x},\label{newv}
\ee
for $x\in (0,1)$, subject to $v^{i+1}(\cdot,0)=v_0,~(v^{i+1})_x(0,t)=0,~v^{i+1}(1,t)=\rho_1$, and

\be
\frac{\p w^{i+1}}{\p t} = \frac{1}{(1-s^i)^2}\frac{\p}{\p x}\left(D(w^i)\frac{\p w^{i+1}}{\p x}\right) + \dot{s}^i\frac{(1-x)}{(1-s^i)}\frac{\p w^{i+1}}{\p x},\label{neww}
\ee
for $x\in (0,1)$, subject to $w^{i+1}(\cdot,0)=w_0,~w^{i+1}(0,t)=\rho_2,~(w^{i+1})_x(1,t)=0$.

These equations are parabolic, by the choice of $(v^i,w^i)$, and therefore have unique smooth solutions $v^{i+1}, w^{i+1}$ on $Q_T$, by Lemma A.1 (see Appendix A), which is a convenient summary of relevant results from classical parabolic theory \cite{lady}; the updated approximants also have corresponding positive diffusivities, by the Maximum Principle.

The coefficients of (\ref{newv}), (\ref{neww}), when expanded into the standard form (\ref{linear_para}), are such terms as $x\dot{s}^i/s^i$, $D(v^i)/(s^i)^2$, $D'(v^i)v^i_x/(s^i)^2$, and analogues for $w^i$, and are therefore dominated in $H^{l,l/2}(Q_T)$ by $|v^i|^{(l+1)}_{Q_T}$ and $|w^i|^{(l+1)}_{Q_T}$. Thus, by Lemma A.1,

\be
|v^{i+1}|^{(l+2)}_{Q_T} \leq  F(C)\left(|v(\cdot,0)|^{(l+2)}_{\Omega}+\rho_1\right),\label{vnew_est1}
\ee
and
\be
|w^{i+1}|^{(l+2)}_{Q_T}\leq F(C)\left(|w(\cdot,0)|^{(l+2)}_{\Omega}+\rho_2\right),\label{wnew_est1}
\ee
where $F$ is some positive function, $\Omega=(0,1)$, and $|\cdot|^{(l+2)}_{\Omega}$ is the ordinary H\"{o}lder norm.

Next, applying Lemma B.1 to (\ref{vnew_est1}) and (\ref{wnew_est1}) results in

\be
|v^{i+1}|^{(l+1)}_{Q_T} - K(\Omega)\|v_0\|_{C^2} \leq K(\Omega)T^{\delta}F(C)\left(|v_0|^{(l+2)}_{\Omega}+\rho_1\right),\label{vnew_est2}
\ee
and
\be
|w^{i+1}|^{(l+1)}_{Q_T} - K(\Omega)\|w_0\|_{C^2} \leq K(\Omega)T^{\delta}F(C)\left(|w_0|^{(l+2)}_{\Omega}+\rho_2\right),\label{wnew_est2}
\ee
where $\delta=\min\left\{\frac{l}{2},\frac{1}{2}(1-l)\right\}$

Thus, if we choose $C>K(\Omega)\max(\|v(\cdot,0)\|_{C^2}, \|w(\cdot,0)\|_{C^2})$, and $T(C)$ is taken sufficiently small, then

\be
\max\left(|v^{i+1}|^{(l+1)}_{Q_T}, |w^{i+1}|^{(l+1)}_{Q_T}\right)\leq C.
\ee

Iteratively, we therefore have that (\ref{vw_c_bound}), (\ref{sdot_bound}), (\ref{vnew_est1}), and (\ref{wnew_est1}) hold uniformly for all $i$.

Next, taking the difference of (\ref{newv}) and the corresponding equation for $v^i$ gives an equation of the form

\bea
\frac{\p}{\p t}(v^{i+1}-v^i) & = & D_i(x,t)\frac{\p^2}{\p x^2}(v^{i+1}-v^i) + a_i(x,t)\frac{\p}{\p x}(v^{i+1}-v^i) + b_i(x,t)\frac{\p}{\p x}(v^i-v^{i-1})\nonumber\\
& + & c_i(x,t)(v^i-v^{i-1}) + d_i(x,t)(s^i-s^{i-1}) + e_i(x,t)(\dot{s}^i-\dot{s}^{i-1}),\label{v_diff_i}
\eea
where $D_i\geq\epsilon$ for all $i$, some $\epsilon>0$, and all the coefficients are bounded in $H^{l, l/2}_{Q_T}$, uniformly in $i$, by the  results just obtained.

For the difference $s^i-s^{i-1}$, we get, from (\ref{rescaled_RH}),

\be
\frac{d}{dt}(s^i-s^{i-1}) = p_i(t)(s^i-s^{i-1}) + q_i(t)(v^i_x-v^{i-1}_x)(1,t) + r_i(t)(w^i_x-w^{i-1}_x)(0,t ),\label{sdiff}
\ee
where $p_i, q_i, r_i$ are bounded in $C^{l/2}$, uniformly in $i$, and from this it is easy to deduce that
\be
|\dot{s}^i-\dot{s}^{i-1}|^{l/2}_{[0,T]}\leq C_1\left(|v^i-v^{i-1}|^{(l+1)}_{Q_T} + |w^i-w^{i-1}|^{(l+1)}_{Q_T}\right)\label{sdiff2}.
\ee

Thus, applying (\ref{sdiff2}) and Lemma A.1 to (\ref{v_diff_i}), we get
\be
|v^{i+1}-v^i|^{(l+2)}_{Q_T}\leq C_1\left(|v^i-v^{i-1}|^{(l+1)}_{Q_T} + |w^i-w^{i-1}|^{(l+1)}_{Q_T}\right).\label{vdiff2}
\ee

Of course, there is also an analogous estimate for $w^{i+1}-w^i$, which, together with (\ref{vdiff2}), implies
\be
|v^{i+1}-v^i|^{(l+2)}_{Q_T} + |w^{i+1}-w^i|^{(l+2)}_{Q_T} \leq C_1\left(|v^i-v^{i-1}|^{(l+1)}_{Q_T} + |w^i-w^{i-1}|^{(l+1)}_{Q_T}\right).
\ee

Finally, an application of Lemma B.1 to the right-hand side of this inequality gives, since all iterates have the same initial data,
\be
|v^{i+1}-v^i|^{(l+2)}_{Q_T} + |w^{i+1}-w^i|^{(l+2)}_{Q_T} \leq C_1T^{\delta}\left(|v^i-v^{i-1}|^{(l+2)}_{Q_T} + |w^i-w^{i-1}|^{(l+2)}_{Q_T}\right).
\ee
(Note that $C_1$ stands for various constants in the above).

If $T$ is chosen so small that $c_1T^{\delta}<1$, it is easy to see that $v^i$ and $w^i$ are Cauchy sequences with respect to  $|\cdot|^{(l+2)}_{Q_T}$-norm. Moreover, it follows from (\ref{sdiff2}) that $s^i$ is Cauchy in $C^{1+\frac{l}{2}}$-norm. By the strong convergence of these sequences, their limits $v, w, s$ satisfy (\ref{v}), (\ref{w}), (\ref{rescaled_RH}) pointwise, and thus constitute a classical solution of $SP_1$.

Turning to the question of uniqueness, suppose we have two smooth solution triples $(v, w, s)$, $(\bar{v}, \bar{w}, \bar{s})$, satisfying (\ref{v}), (\ref{w}), (\ref{rescaled_RH}), and having the same initial data. Subtracting the equation satisfied by $\bar{v}$ from that satisfied by $v$ then gives an equation of the form

\be
\p_t(v-\bar{v}) = D(x,t)\p_x^2(v-\bar{v}) + a(x,t)\p_x(v-\bar{v}) + b(x,t)(v-\bar{v}) + c(x,t)(s-\bar{s}) + d(x,t)(\dot{s}-\dot{\bar{s}}),\label{v-vbar} 
\ee
where the coefficients are smooth, and $D(x,t)\geq\epsilon$, some $\epsilon>0$. An analogous equation is obtained for $w-\bar{w}$ in the same way.

Thus, by Lemma A.1, $v-\bar{v}$ can be bounded by the inhomogeneity in (\ref{v-vbar}), according to

\bea
|v-\bar{v}|^{(l+2)}_{Q_T} & \leq & C\left(|c(x,t)(s-\bar{s})|^{(l)}_{Q_T} + |d(x,t)(\dot{s}-\dot{\bar{s}})|^{(l)}_{Q_T} \right)\nonumber\\
& \leq & C\left(|s-\bar{s}|^{(l)}_{Q_T} + |\dot{s}-\dot{\bar{s}}|^{(l)}_{Q_T}\right)\nonumber\\
& \leq & C\left(|v_x-\bar{v}_x|^{(l)}_{Q_T} + |w_x-\bar{w}_x|^{(l)}_{Q_T}\right)\nonumber\\
& \leq & C\left(|v-\bar{v}|^{(l+1)}_{Q_T} + |w-\bar{w}|^{(l+1)}_{Q_T}\right),\label{inequality4}
\eea
for various constants $C$, where we used the Rankine-Hugoniot condition (\ref{rescaled_RH}) to get the third line.

Adding the analogous inequality for $w-\bar{w}$ therefore results in
\be
|v-\bar{v}|^{(l+2)}_{Q_T} + |w-\bar{w}|^{(l+2)}_{Q_T} \leq  C\left(|v-\bar{v}|^{(l+1)}_{Q_T} + |w-\bar{w}|^{(l+1)}_{Q_T}\right),
\ee
which, with the help of Lemma B.1, implies

\be
|v-\bar{v}|^{(l+2)}_{Q_T} + |w-\bar{w}|^{(l+2)}_{Q_T} \leq  CT^{\delta}\left(|v-\bar{v}|^{(l+2)}_{Q_T} + |w-\bar{w}|^{(l+2)}_{Q_T}\right).
\ee 
Choosing $T$ small enough gives, finally, $v(x,t)=\bar{v}(x,t),~w(x,t)=\bar{w}(x,t)$, as required $\square$

\subsubsection{Global existence subject to a sign condition on the solution gradient at the discontinuity}

We will prove a global-existence theorem for $SP_1$ by showing that, for a local classical solution $(\rho(x,t),s(t))$, the gradient, $\rho_x(x,t)$, and hence also $\dot{s}(t)$, are {\em a priori} bounded, provided that the one-sided limits $\rho_x(s^{\pm}(t),t)$ satisfy a sign condition. The proof is rather different from, and somewhat shorter than, that presented in \cite{meir}.

First, with the change of dependent variable $\sigma=K(\rho)$, which is smooth and invertible in each phase, the governing PDE becomes

\be
\sigma_t = \widetilde{D}(\sigma)\sigma_{xx},\label{global_sigma}
\ee
where $\widetilde{D}(\sigma)=D(K^{-1}(\sigma))>0$, the boundary conditions become 
\bd
\sigma_x(0)=\sigma_x(1)=0,\quad\sigma(s^-)=K(\rho_1):=\sigma_1,\quad\sigma(s^+)=K(\rho_2):=\sigma_2,
\ed
and the Rankine-Hugoniot condition takes the form
\be
\frac{ds}{dt} = -\frac{(\sigma_x(s^+) - \sigma_x(s^-))}{(\rho_2-\rho_1)}.\label{RH_sigma_0}
\ee

Next, note that integration by parts gives, for any smooth $\sigma$, and $m=1,2,\ldots$,

\be
\int_0^{s(t)} \p_t(\sigma_x)^{2m}~dx = -\int_0^{s(t)} 2m(m-1)\sigma_x^{2m-2}\sigma_{xx}\sigma_t~dx + \left[2m\sigma_x^{2m-1}\sigma_t\right]_0^{s(t)}\label{int_by_parts}.
\ee

Assuming now that $(\sigma,s)$ solves $SP_1$, we have $\sigma(s^-(t),t)=\sigma_1$, and hence, by differentiation,\linebreak $\sigma_x\dot{s} + \sigma_t=0$ at $x=s^-$, which, with the aid of (\ref{int_by_parts}), leads to

\bea
\frac{d}{dt}\|(\sigma_x)^m\|^2_{L^2(0,s)} & = & -\int_0^s 2m(2m-1)(\sigma_x)^{2m-2}\sigma_{xx}\sigma_t~dx + (1-2m)(\sigma_x)^{2m}\dot{s}\left.\right|_{x=s^-}\nonumber\\
 & \leq & (1-2m)\dot{s}(t)(\sigma_x)^{2m}|_{x=s^-},\label{sig_left}
\eea
where we used (\ref{global_sigma}) to discard the integral term.

By a similar calculation, there also follows

\be
\frac{d}{dt}\|(\sigma_x)^m\|^2_{L^2(s,1)} \leq -(1-2m)\dot{s}(t)(\sigma_x)^{2m}|_{x=s^+}.\label{sig_right}
\ee

If we regard $\sigma_x$ as a function in $L^2(0,1)$, then adding (\ref{sig_left}) and (\ref{sig_right}), and using (\ref{RH_sigma_0}), results in

\be
\frac{d}{dt}\|(\sigma_x)^m\|_2^2 \leq (1-2m)\frac{(\sigma_x(s^+)-\sigma_x(s^-))}{(\rho_2-\rho_1)}((\sigma_x)^{2m}(s^+)-(\sigma_x)^{2m}(s^-)).\label{sig_decay}
\ee

Now, if $2m=2^n,~n=1,2,\ldots$ then, inductively, there is a {\em positive} multinomial $F_n(a,b)$ such that 

\be
(a-b)(a^{2^n}-b^{2^n}) = F_n(a,b)(a+b),\quad\forall a,b.\label{F_n}
\ee  

Using (\ref{F_n}) in the right-hand side of (\ref{sig_decay}) with $a=\sigma_x(s^+), b=\sigma_x(s^-)$ shows that $\|(\sigma_x)^m\|_2^2$ is decreasing for $2m=2^n$, provided 
\be
\sigma_x(s^+) + \sigma_x(s^-)\geq 0\quad\forall t.\label{grad_sign}
\ee

This condition can be guaranteed, for example, by choosing data $(\rho_0, s_0)$ for $SP_1$ satisfying the traditional condition (\ref{class_data}), and by appealing to the Maximum Principle in each phase.

Thus, assuming (\ref{grad_sign}) holds, we have

\be
\|\sigma_x\|_{L^{2^n}}(t) \leq \|\sigma_x\|_{L^{2^n}}(0) \leq \|\sigma_x\|_{\infty}(0),
\ee
for $n=1,2,\ldots$, and therefore, by Theorem 2.8 of \cite{adams},

\be
\|\sigma_x\|_{\infty}(t)\leq \|\sigma_x\|_{\infty}(0).
\ee

In 1-d, such an {\em a priori} bound on $\|\nabla\sigma\|_{\infty}$, which of course also gives a pointwise bound on the corresponding $\rho_x$ and on $\dot{s}(t)$, is actually enough to continue the local solution of $SP_1$ obtained in Sect.2.3, by standard theory. Indeed, if $v$ and $w$ are as in (\ref{v})-(\ref{w}), then combining the estimates of Appendix A gives, for $3<q<4$, and a constant $C$ which is controlled by $\nabla\sigma$,
\be
|v|^{\left(2-\frac{3}{q}\right)}_{Q_T}\leq C\left(\|v\|_{C^2((0,1))}+\rho_1\right),
\ee
and an analogous estimate for $w$.

Thus, for some $\beta>0$, the $H^{\beta+1,\frac{1}{2}(\beta+1)}$-norm of $v$ and $w$ is controlled by $\nabla\sigma$ and the $C^2$-norm of the initial data. Given this, Lemma A.1 implies in turn that the $H^{\beta+2,\frac{1}{2}(\beta+2)}$-norm of $v$ and $w$ is controlled by $\nabla\sigma$ and the $C^{\beta+2}$-norm of the initial data. Thus, the local classical solution can always be extended onto to a larger time interval, provided, of course, $s(t)$ stays away from the domain boundary.

We have already seen (Theorem 2.1) that if the total mass $M$ lies outside the interval $(\rho^{\flat},\rho^{\sharp})$, it cannot be ruled out that $s(t)$ will hit the domain boundary in finite time; should this occur, the solution can be continued via the Neumann problem on $\Omega = (0,1)$, as mentioned above. Note that, on a disappearing phase, the solution (in terms of the original variable $\rho$) merely approaches $\rho_1$ or $\rho_2$, as appropriate, by the boundedness of the gradient. Since we have already extended our definition of $SP_1$ to cover such eventualities, we have therefore proved

\begin{theorem}
Given phase-wise smooth initial data $(\rho_0(x),s_0)$ satisfying (\ref{class_data}) and the first-order compatibility condition, problem $SP_1$ has a unique smooth, global-in-time solution $(\rho(x,t),s(t))$, such that the corresponding $\|\sigma_x\|_{\infty}$ is monotonically decreasing for all $t$.
\end{theorem}

In general, it is not clear whether a corresponding global-existence result can be obtained in the situation where (\ref{class_data}) is not satisfied. An exception to this is the special case of 1-phase problems, to which we now turn our attention.

\subsubsection{Global existence for 1-phase problems; Mass Lagrange coordinates}

If the initial density is constant for either $x<s_0$ or $x>s_0$, then we refer to the corresponding evolution problem $SP_1$ as a one-phase problem. As a consequence of the Neumann condition at the domain boundary, it then turns out that there is a change of the spatial variable, different from a simple rescaling, which transforms the one-phase $SP_1$ into a regular quasilinear parabolic problem on a fixed domain. This allows us to prove a global-existence theorem, using standard parabolic theory, regardless of the direction of the gradient at the phase boundary. 

Without loss of generality, let us assume that $\rho$ is constant ($=\rho_2$) in the high-density phase, and variable in the low-density phase, and let us again make the change of dependent variable $\sigma=K(\rho)$ on $0\leq x\leq s(t)$. Since we are assuming $\rho<\rho^{\flat}$ in the low-density phase, this change of variable is invertible there, with inverse denoted by $\rho=b(\sigma):=K^{-1}(\sigma)$.

In terms of $\sigma$, the governing equation (\ref{cont_rho}) becomes
\be
\frac{\p b(\sigma)}{\p t} = \Delta\sigma,\label{b_sigma}
\ee
subject to $\sigma_x=0$ at $x=0$, and $\sigma=K(\rho_1)$ at $x=s(t)$, while the Rankine-Hugoniot condition takes the form
\be
\frac{ds}{dt}=\frac{\p\sigma}{\p x}(s^-)/(\rho_2-\rho_1).
\ee

Next, we transform (\ref{b_sigma}) by introducing so-called Mass Lagrange coordinates \cite{meir}, $(\tau,y)$, which are defined by

\be
\tau=t,\qquad y=\int_x^{s(t)}\limits\left[\rho_2-b(\sigma(\hat{x},t))\right]~d\hat{x}.
\ee
Since, by construction, $\frac{\p y}{\p x}=-(\rho_2-b(\sigma))<0$, this coordinate change is good, and, introducing $v(y,\tau)=\sigma(x,t)$, (\ref{b_sigma}) now takes the form
\be
\frac{b'(v)}{(\rho_2-b(v))^2}\frac{\p v}{\p\tau}=\frac{\p^2v}{\p y^2},\label{b_v}
\ee
which is a regular quasilinear parabolic equation for $v$.

Clearly, the moving boundary $x=s(t)$ gets mapped to $y=0$, and, as a consequence of conservation of mass, the domain boundary $x=0$ gets mapped to 
\be
y_0:=\int_0^{s(0)}\limits[\rho_2-b(\sigma(x,0))]~dx=\textrm{const.}.
\ee
Thus, (\ref{b_v}) is to be solved on a fixed spatial domain subject to the boundary conditions $\frac{\p v}{\p y}=0$ at $y=y_0$, and $v=K(\rho_1)$ at $y=0$.

It is well known that this problem has a unique, global classical solution $v(y,\tau)$, given $H^1$ initial data compatible with the Dirichlet condition at $y=0$ \cite{amann,lady}. The corresponding jump location $s(t)$, which can be reconstructed from $v(y,\tau)$ by integrating, for each $t=\tau$, the equation
\be
\frac{\p y}{\p x}=b(v(y,\tau))
\ee
from the point $(x=0,~y=y_0)$ until $y$ hits zero, could conceivably hit $x=1$ in finite time (see the numerics in the next subsection). If this occurs, the solution of (\ref{b_v}) should be stopped, and then continued for all time via NP for (\ref{cont_rho}) on $\Omega=(0,1)$.

We thus arrive at

\begin{theorem}
Given smooth, one-phase initial data which satisfies the first-order compatibility condition, problem $SP_1$ has a unique smooth, global-in-time solution $(\rho(x,t),s(t))$.
\end{theorem}

\subsubsection{Linear stability near a discontinuous steady state}

The use of rescaled coordinates, as introduced early on in this subsection, allows us to investigate linear stability of the unique discontinuous steady state (\ref{s_star}) of $SP_1$, which exists as long as \linebreak $\rho_1<M<\rho_2$.

Linearising (\ref{v}) and (\ref{w}) around $v=\rho_1$ and $w=\rho_2$, respectively, and using hats to denote differentials, we get the pair of heat equations
\bea
\hat{v}_t & = & \frac{1}{(s^{\ast})^2}D(\rho_1)\hat{v}_{xx},\\
\hat{w}_t & = & \frac{1}{(1-s^{\ast})^2}D(\rho_2)\hat{w}_{xx},
\eea
which are to be solved subject to $\hat{v}(1)=0$, $\hat{v}_x(0)=0$, $\hat{w}(0)=0$, $\hat{w}_x(1)=0$, while linearising (\ref{resc_COM}) around $s=s^{\ast}$ gives
\be
\hat{s}=(\rho_1-\rho_2)^{-2}\left\{-(\rho_1-\rho_2)\int_0^1\hat{w}~dx + (\rho_2-M)\int_0^1(\hat{v}-\hat{w})~dx\right\}.
\ee

For a solution $(\hat{v},\hat{w},\hat{s})$ of the linearisation, it is therefore clear that, say, $\|\hat{v}\|_{\infty}\rightarrow 0$ and  $\|\hat{w}\|_{\infty}\rightarrow 0$ as $t\rightarrow\infty$, and also that $\hat{s}\rightarrow 0$ as $t\rightarrow\infty$. In this sense, then, the unique discontinuous steady state of $SP_1$ is always linearly stable, whenever it exists.

\subsection{Numerical simulations}

As well as facilitating mathematical analysis, the rescaled, fixed-boundary representation of $SP_1$, (\ref{v}), (\ref{w}), (\ref{rescaled_RH}), also comes in useful for numerical simulations. 

Specifically, we use a method-of-lines approach in which the diffusion terms in (\ref{v}) and (\ref{w}) are discretised using the random-walk model (\ref{walk}), while the advection terms are discretised by means of a standard, explicit upwinding scheme. A simulation in which the global-existence criterion (\ref{class_data}) is satisfied is shown in Figure 2, and one in which it is violated is shown in Figure 3. In each case, the solution approaches the appropriate discontinuous steady state at large times. Several other simulations have been carried out in the case where (\ref{class_data}) is violated, and no singularities have been observed to develop.

\begin{figure}

\centering

\resizebox{5in}{5in}{\includegraphics{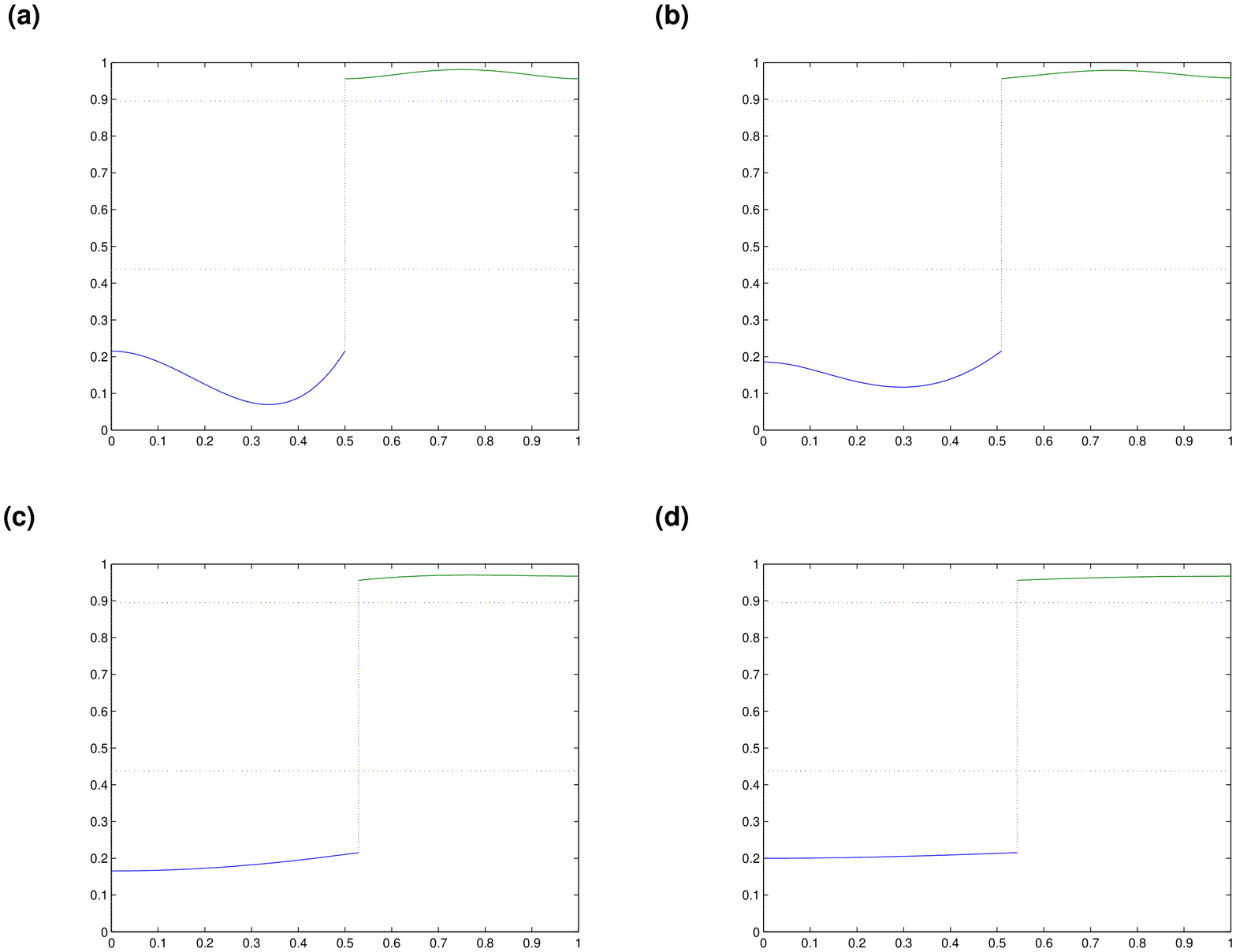}}

\caption{Solution of $SP_1$ with $\alpha=0.85$ and initial data satisfying (\ref{class_data}). Density profile shown at (a) t=0, (b) t=0.0125, (c) t=0.1289, (d) t=0.456.}

\end{figure}

\begin{figure}

\centering

\resizebox{5in}{5in}{\includegraphics{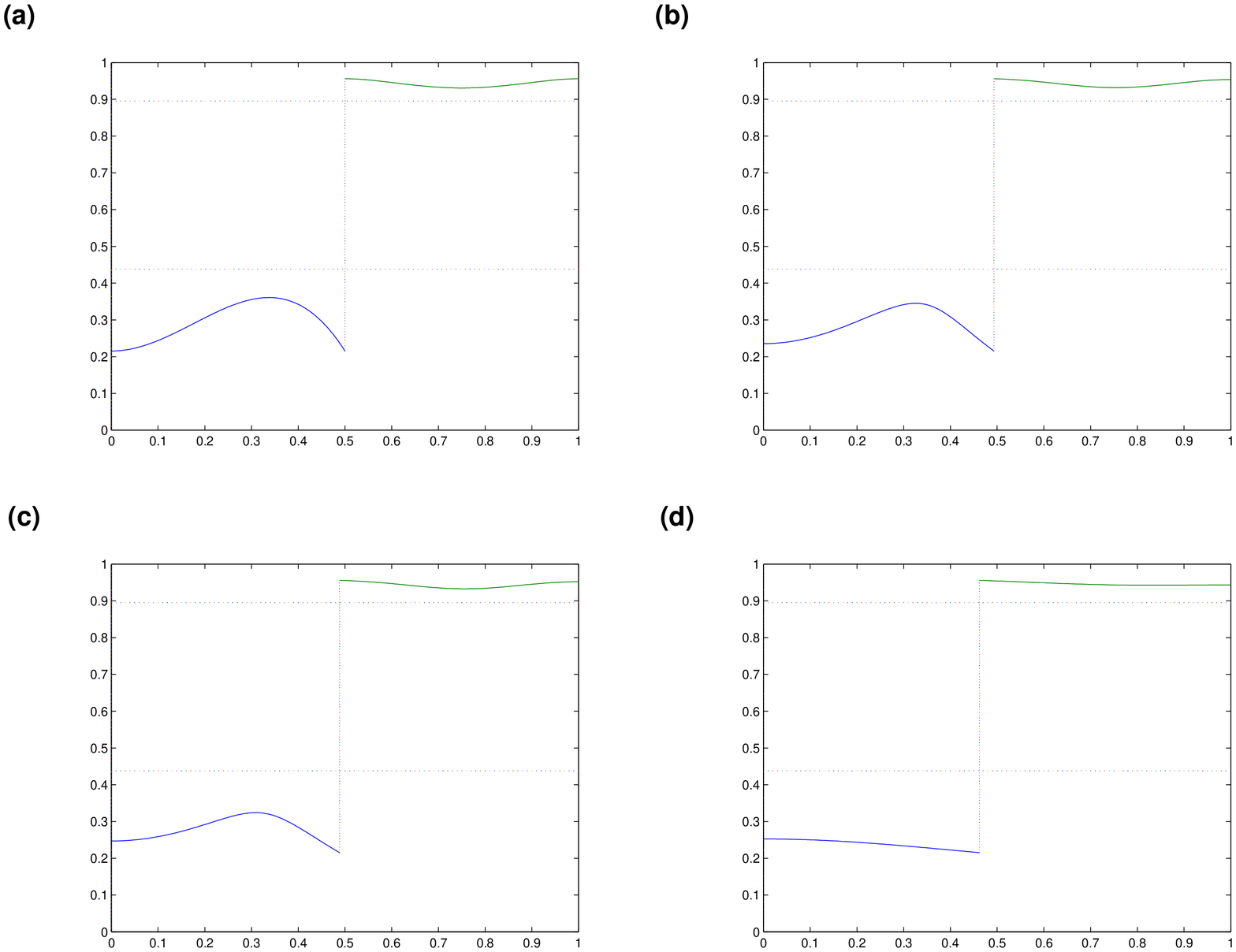}}

\caption{Solution of $SP_1$ with $\alpha=0.85$ and initial data violating (\ref{class_data}). Density profile shown at (a) t=0, (b) t=0.0125, (c) t=0.1289, (d) t=0.456.}

\end{figure}

Turning to the question of bistability for cetain values of the mass $M$, we next show in Figure 4 a close-up of a solution with $M=0.3184$ and $\alpha=0.85$, for which the phase boundary hits the domain boundary in finite time. This should be compared with the simulation shown in Figure 5, in which the solution has the same mass and adhesion coefficient, but this time evolves to the unique discontinuous steady state as $t\rightarrow\infty$. 

It is perhaps worth noting that such bistability cannot occur if the traditional restriction (\ref{class_data}) is imposed on the initial data, by Theorem 2.1.

\begin{figure}

\centering

\resizebox{5in}{5in}{\includegraphics{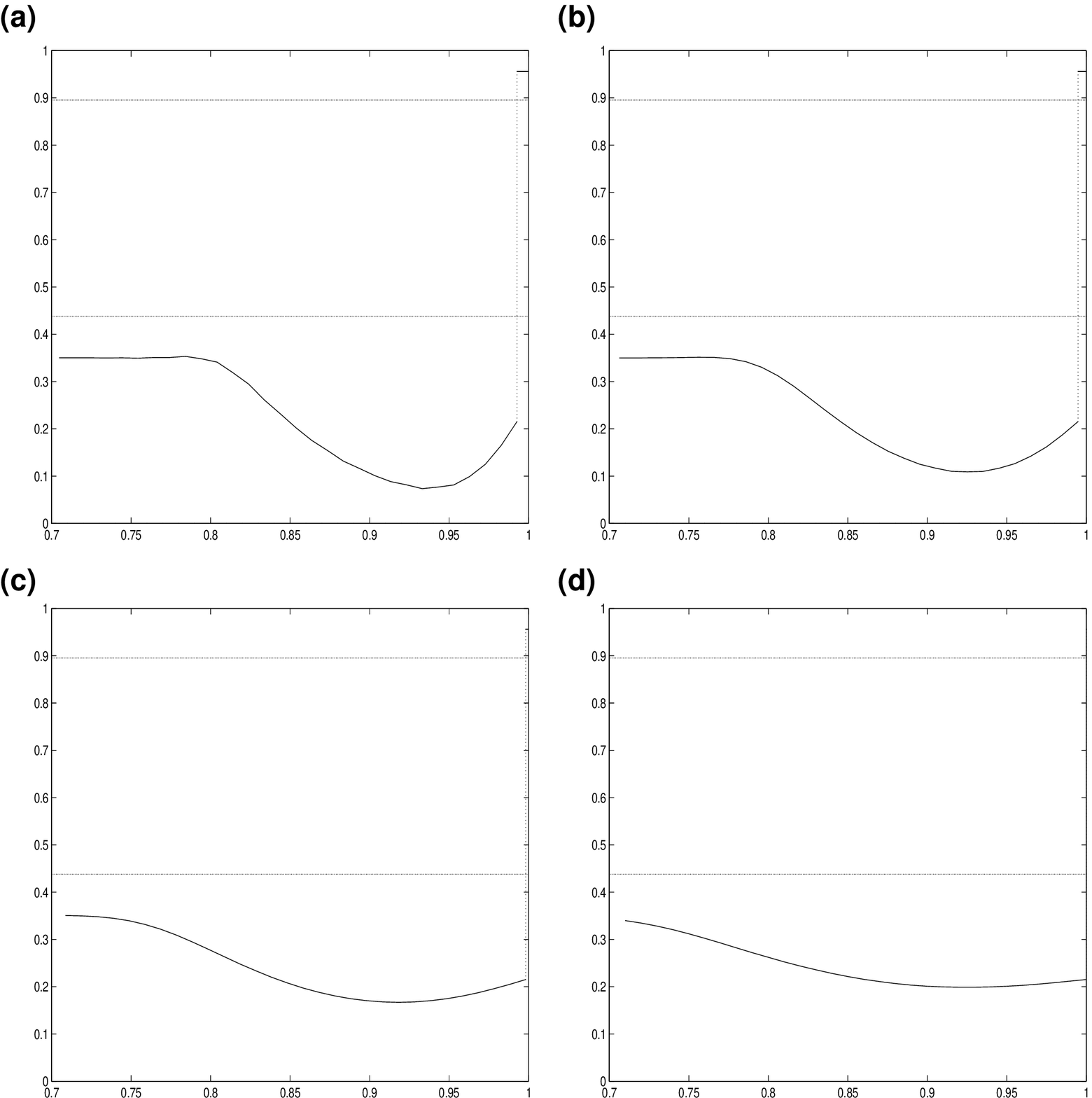}}

\caption{Close-up of a solution of $SP_1$ with $M=0.3184$ and $\alpha=0.85$, such that $s(t)$ hits $x=1$ in finite time. Density profile shown at (a) t=0, (b) t=0.001, (c) t=0.0051, (d) t=0.012. The density to the left of $x=0.7$ is essentially constant throughout the simulation.}

\end{figure}

\begin{figure}

\centering

\resizebox{5in}{5in}{\includegraphics{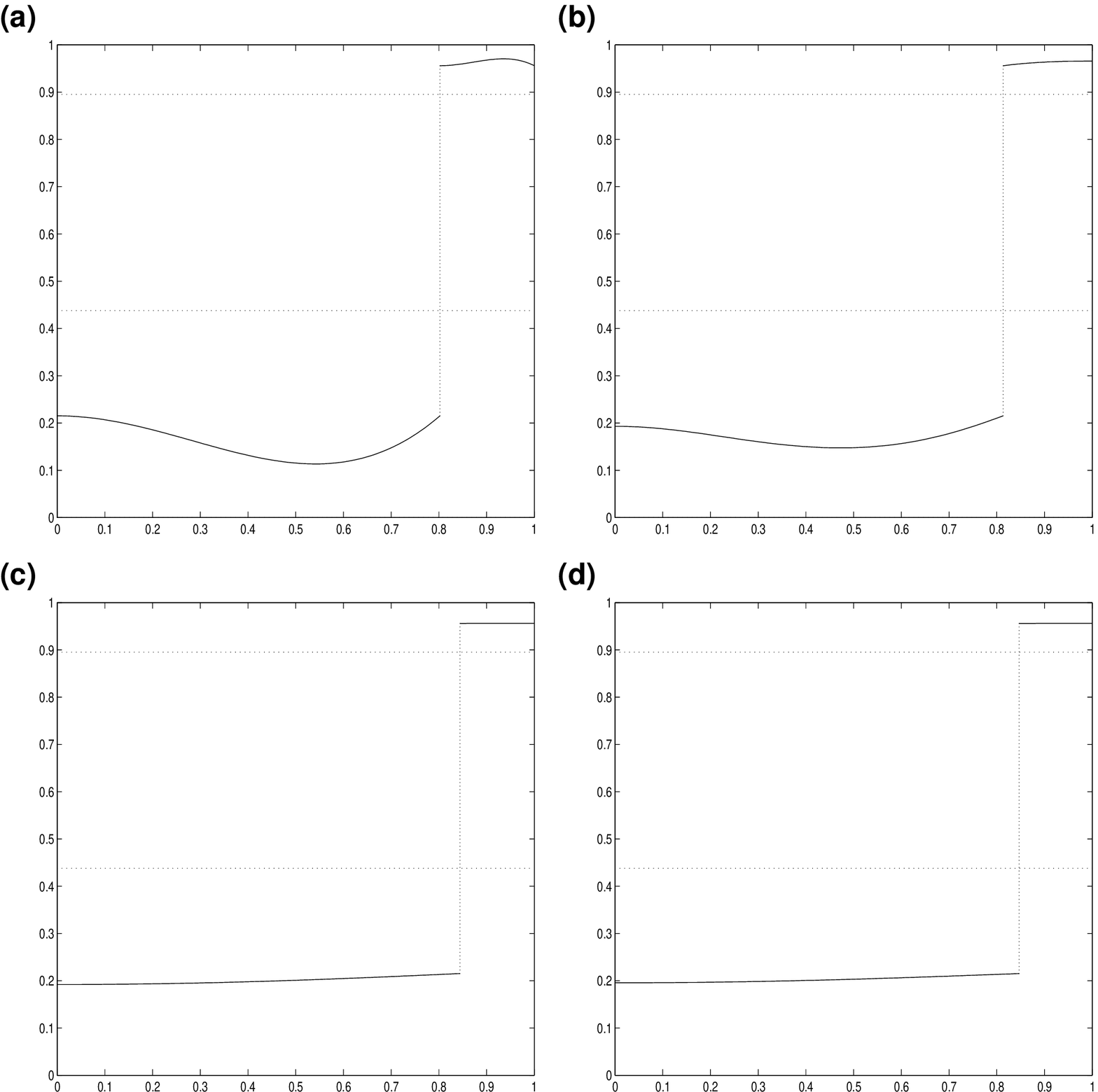}}

\caption{Solution of $SP_1$ with $M=0.3184$ and $\alpha=0.85$ which evolves to a discontinuous steady state as $t\rightarrow\infty$. Density profile shown at (a) t=0, (b) t=0.0373, (c) t=0.6264, (d) t=0.7479.}

\end{figure}

\section{Multi-phase problems}

We now consider the general case of our Stefan-problem set-up, in which the initial density profile jumps $n$ times across the the unstable region $I_{\alpha}$. This initial/boundary-value problem will be denoted by $SP_n$.

\subsection{Formulation}

The problem $SP_n$ consists of looking for a piecewise-smooth $\rho$ which satisfies the diffusion equation (\ref{cont_rho}) away from $n$ discontinuities $s_i(t)$, with $0<s_i<s_{i+1}<1$, such that $\rho$ always jumps between $\rho_1$ and $\rho_2$ at the $s_i$, which evolve according to the Rankine-Hugoniot condition
\be
\frac{ds_i}{dt} = \frac{(J_i^+ - J_i^-)}{[\rho_i]},\label{RH_i}
\ee
where $J_i^{\pm}=-D(\rho(s_i^{\pm}(t),t))\p_x\rho(s_i^{\pm}(t),t)$, and $[\rho_i]=(\rho(s_i^+,t)-\rho(s_i^-,t))$~($=\pm(\rho_2-\rho_1)$), and the Neumann condition is again imposed at the domain boundary $x=0,1$; a schematic for $SP_n$ is given in Figure 6.

We also allow for the possibility that $s_1$ or $s_n$ may hit the domain boundary in finite time, or that neighbouring discontinuities could collide, leading to the annihilation of a phase. Should any of these events occur, the solution can be continued via $SP_{n-1}$ or $SP_{n-2}$, as appropriate, and so on. A solution continued in this way will still be referred to globally as a solution of $SP_n$.

\subsection{Steady states and their stability}

The steady-state picture for $SP_n$, given $n\geq 2,$ is more complicated than that for $SP_1$. 

First of all, for a given mass $M$ satisfying $\rho_1<M<\rho_2$, there is a continuum of two-valued, $n$-jump steady states, each of which is given by a choice of the $s_i$ which merely has to be compatible with $M$. Moreover, discontinuous steady states with fewer than $n$ jumps can also be considered permissible - these could be approached dynamically by (multiple) coalescence events, and/or by (successive) collisions of phase boundaries with the domain boundary. Each steady state with more than one discontinuity is expected to be merely neutrally stable, since the total mass is invariant under small translations of an internal phase (i.e., one which does {\em not} touch the domain boundary).

If the stronger condition $\rho^{\flat}<M<\rho^{\sharp}$ holds, then there is no possible continuous steady state, and thus, for a global solution of $SP_n$, at least one discontinuity must remain as $t\rightarrow\infty$. If, instead, $\rho_1<M\leq\rho^{\flat}$ or $\rho^{\sharp}\leq M<\rho_2$, then the uniform steady state exists alongside the discontinuous family already discussed. In this regard, note that Figures 4 and 5 could be considered as simulations of one half of a reflection-symmetric $SP_2$, with $\rho_1<M<\rho^{\flat}$, in which the central high-density phase is either annihilated in finite time (Figure 4) or preserved as $t\rightarrow\infty$ (Figure 5), depending on the proximity of the initial data to the (unique) discontinuous, symmetric steady state.

In general, it is difficult to say anything analytical about the stability of steady states of $SP_n$, although in each phase the solution will still decay to $\rho_1$ or $\rho_2$, as appropriate, exponentially in $L^2$, by essentially the same calculation as in Section 2.2, as long as the relevant phase boundaries stay away from the domain boundary.

Finally, note that if $M\geq\rho_2$ or $M\leq\rho_1$, then there is only one possible steady state, namely the uniform one, and, by $L^2$-decay, phase boundaries for a global solution must disappear in finite time (or possibly infinite time in the exceptional cases $M=\rho_1, \rho_2$) via coalescence events, or by merging with the domain boundary. For the subsequent $NP$, exponential $L^2$-convergence to the uniform steady state $\rho=M$ follows as before.

\begin{figure}

\centering

\resizebox{5in}{5in}{\includegraphics{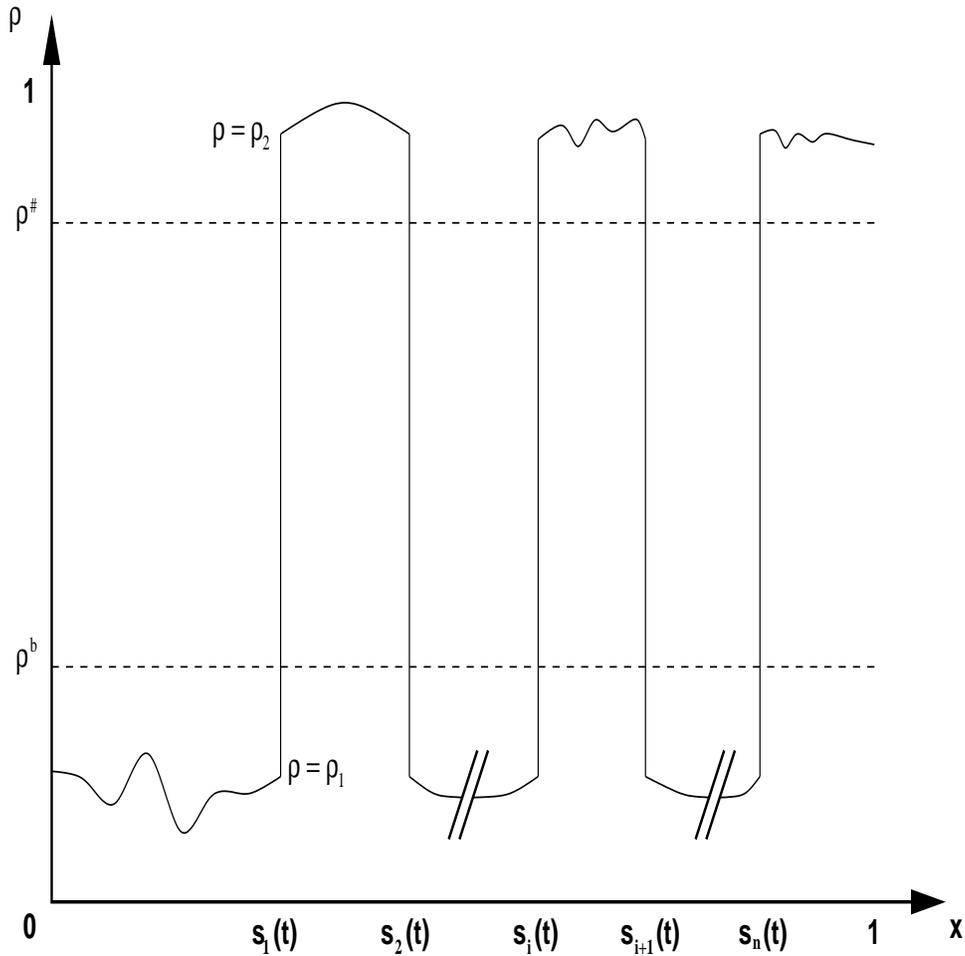}}

\caption{Example of a density profile for the multi-jump Stefan problem $SP_n$.}

\end{figure}

\subsection{A local-existence-and-uniqueness theorem for smooth data}

The multi-phase problem, $SP_n$, is solved, locally in time, in the same way as $SP_1$; in each phase the spatial variable is rescaled in order to fix the moving boundary (or boundaries), and the same estimates go through as before.

To be explicit, on the interval $[s_i,s_{i+1}]$ we make the coordinate transformation
\be
(x,t)\rightarrow(\hat{x},\hat{t}) :\quad\hat{x}=\frac{(x-s_i)}{(s_{i+1}-s_i)},\quad\hat{t}=t,
\ee
such that, writing $v_i(\hat{x},\hat{t})=\rho(x,t)$, the governing PDE (\ref{cont_rho}) becomes, upon dropping hats
\be
\frac{\p v_i}{\p t} = (s_{i+1}-s_i)^{-2}\frac{\p}{\p x}\left(D(v_i)\frac{\p v_i}{\p x}\right) + \frac{\left((\dot{s}_{i+1}-\dot{s}_i)x+\dot{s_i}\right)}{(s_{i+1}-s_i)}\frac{\p v_i}{\p x},\label{v_i}
\ee
on $[0,1]\times[0,T]$, for $i=0,1,2,\ldots,n$, where we adopt the conventions $s_0(t)=0$ and $s_{n+1}(t)=1$.

On an internal phase, this equation is to be solved subject to the appropriate Dirichlet conditions at $x=0,1$, while on an extremal phase the mixed Neumann/Dirichlet conditions are used, as in the two-phase case (see Figure 6).

The rescaled family of Rankine-Hugoniot conditions takes the form
\be
\frac{ds_i}{dt} = -\left(\frac{D(\rho_2)(v_{i+1})_x(0,t)}{(s_{i+1}-s_i)}-\frac{D(\rho_1)(v_i)_x(1,t)}{(s_i-s_{i-1})}\right)(\rho_2-\rho_1)^{-1},\label{rescaled_RH_i}
\ee
for $i=1,2,\ldots,n$.

Equations (\ref{v_i}), (\ref{rescaled_RH_i}) are solved by the same kind of iteration employed for $SP_1$; the required H\"{o}lder estimates are obtained via linear parabolic theory and Picard's Theorem for systems of ODEs.

The upshot of all this is

\begin{theorem}
Given initial data $(v_i^0, w_i^0, s_i^0)$, $i=1,2,\ldots,n$, for which the $v_i^0$ and $w_i^0$ belong to the H\"{o}lder space $C^{l+2}$, $l>0$, and such that the  first-order compatibility conditions are satisfied, the system of equations (\ref{v_i}), (\ref{rescaled_RH_i}), and hence problem $SP_n$, has a unique classical solution on some small time interval $[0,T]$.
\end{theorem}

\subsection{Continuation of the local solution, subject to a sign condition on the solution gradient at discontinuities}

With the same notation as in Sect. 2.4, and by a similar calculation, it is straightforward to see that the gradient of a local solution of $SP_n$ satisfies, for $m=1,2,\ldots$, the {\em a priori} estimate
\be
\frac{d}{dt}\|(\sigma_x)^m\|^2_{L^2((0,1))}\leq\left(\frac{1-2m}{\rho_2-\rho_1}\right)\sum_i\mathrm{sign}[\rho]_i\left(\sigma_x(s_i^+)-\sigma_x(s_i^-)\right)\left(\sigma_x^{2m}(s_i^+)-\sigma_x^{2m}(s_i^-)\right).
\ee

Thus, $\|\sigma_x\|_{\infty}$, and hence also the $\|\rho_x\|_{L^{\infty}((s_i,s_{i+1}))}$ and $\dot{s}_i$, are bounded for all time, provided
\be
\mathrm{sign}[\rho]_i\left(\sigma_x(s_i^+)+\sigma_x(s_i^-)\right)\geq 0,\quad\forall i.
\ee

This condition holds if, for example,

\be
\rho_0\leq\rho_1\quad\textrm{in low-density phases},\quad\rho_0\geq\rho_2\quad\textrm{in high-density phases},\label{grad_cond2}
\ee
which is the multi-phase analogue of (\ref{class_data}).

Inequalities (\ref{grad_cond2}) therefore guarantee global existence of the corresponding solution to $SP_n$, modulo coalescence events, and the possibility that an extremal discontinuity might hit the boundary in finite time. 

\subsection{Coalescence events, and continuation thereafter}

In order to continue a solution of $SP_n$ after the coalescence of two phases, it is necessary to extend our existence theory (which has thus far required the initial data to lie in $C^{l+2}$) to the case where the initial density profile may have a `corner' in one of the phases. This can be seen by considering the situation illustrated in Figure 7, in which two high-density phases coalesce, thus annihilating a low-density phase. 

\begin{figure}

\centering

\resizebox{6in}{5in}{\includegraphics{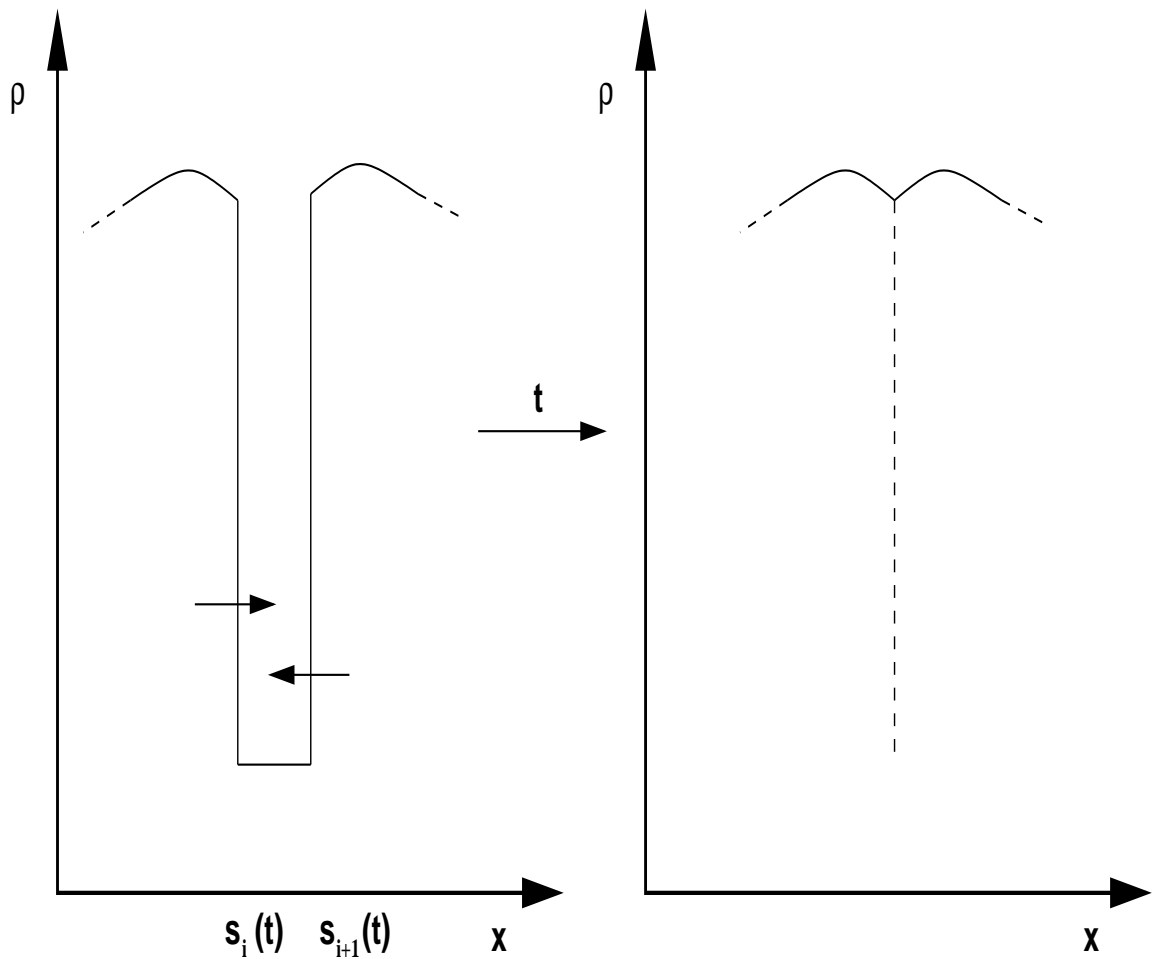}}

\caption{Two coalescing high-density phases; note that the resulting single high-density phase has a corner at the point and time of coalescence.}

\end{figure}

If we make the assumption that condition (\ref{grad_cond2}) holds, then, in fact, the required extension for such phase-wise $H^1$ data follows from Theorem 14, Sect.V of \cite{meir}, and hence, putting together the results of this section, we therefore arrive at
\begin{theorem}
Given $n$ initial jump locations $(s_i)_0$, and a phase-wise smooth initial density profile $\rho_0$ satisfying the first-order compatibility conditions at phase boundaries, along with the gradient condition (\ref{grad_cond2}), the multi-phase Stefan problem $SP_n$ has a unique, global classical solution, on the understanding that phases may in time be annihilated via coalescence events, or that an extremal discontinuity may hit the domain boundary in finite time. Moreover, $\|\sigma_x\|_{\infty}$ is monotonically decreasing for all time. 
\end{theorem}

We end by noting that it is not clear whether even a local existence theorem for $SP_n$ can be proved for merely $H^1$-data if (\ref{grad_cond2}) is not satisfied; one can, for example, proceed by approximating $H^1$ data with smooth data, but it seems that the {\em a priori} bound on $\sigma_x$ implied by (\ref{grad_cond2}) is required to get the necessary convergence.

\section{Concluding remarks}

The results of this paper can be interpreted as saying that cell-cell adhesion is enough to stabilise sharp-edged aggregations of diffusing cells if there is enough mass present in the biological domain, and that diffusion homogenises the cell density in each given high- or low-density phase. For low masses, however, narrow high-density regions can become annihilated even if the adhesion is very strong - that is to say, diffusion in the low-density regions wins, leading to a globally uniform cell density at large times. In the case of intermediate masses, bistability becomes possible - either adhesion or diffusion can win out, depending on the profile of the initial data. 

Next, while being motivated by the adhesion-diffusion equations (\ref{cont_rho}) and (\ref{mod-equ}), it should be clear that the analysis carried out in this paper does not require that the diffusivity $D(\rho)$ have the special form (\ref{D}); indeed, all arguments go through for {\em any} equation of the form $\rho_t=K(\rho)_{xx}$, provided the $C^2$ function $K$ is increasing outside an unstable interval of $\rho$ values, and provided  $[K(\rho)]_i=0$ at jump locations $s_i$.

Finally, one rather obvious biologically-relevant extension of the work described here would be to carry out a similar analysis with a chemotactic term factored into the right-hand side of (\ref{cont_rho}). In that situation, one could imagine beginning with a low-density initial datum, evolving the solution until, through chemotactic aggregation,~$\rho$ hits the unstable region $I_{\alpha}$ at some point $x_c$, and then continuing the solution via a Stefan problem with an initial spike at $x_c$ which jumps from $\rho_1$ to $\rho_2$. The transition from well-posed Neumann problem to Stefan problem is somewhat singular in that case, and is the subject of ongoing analytical investigation. 

\section*{Acknowledgement}

The author wishes to thank Christian Schmeiser for suggesting this problem, and for several helpful discussions along the way.

\section*{Appendix A. Linear parabolic theory}

First of all, we have

\textbf{Lemma A.1}~~~~{\em The equation}

\be
\frac{\p u}{\p t} = D(x,t)\frac{\p^2 u}{\p x^2} + a(x,t)\frac{\p u}{\p x} + b(x,t)u + f(x,t),\label{linear_para}
\ee
{\em with initial data $u_0\in C^{l+2}([0,1])$, $0<l<1$, subject to the boundary conditions $u_x(0,t)=0$, $u(1,t)=u_1$, and the first-order compatibility condition}
\be
D(x,0)\frac{\p^2 u_0}{\p x^2} + a(x,0)\frac{\p u_0}{\p x} + b(x,0)u_0 + f(x,0) = 0,
\ee
{\em has a unique solution $u\in H^{l+2,(l+2)/2}(Q_T)$ on $Q_T=[0,1]\times[0,T]$, which satisfies the estimate}

\be
|u|^{(l+2)}_{Q_T}\leq c\left(|f|^{(l)}_{Q_T} + \|u_0\|_{C^{l+2}([0,1])} + |u_1|\right),
\ee
{\em with the constant $c$ remaining bounded as $T\rightarrow 0$, provided $D(x,t)\geq\epsilon>0$ and all the coefficients and the inhomogeneity in (\ref{linear_para}) are bounded in $H^{l,l/2}(\Omega_T)$.}
\\~\\
For the proof, see \cite{lady}, Ch.4, Thms. 5.2-5.4.

Next, let $H^2((0,1))$ denote the usual second-order $L^2$-type Sobolev space on the unit interval, let $W^{1,q}((0,1))$, $q>1$, be the Sobolev space with norm
\be
\|u\|_{1,q} = \|u\|_{L^q((0,1))} + \|u_x\|_{L^q((0,1))},
\ee

and let $W^{2-\frac{2}{q},q}((0,1))$ denote the fractional-order Sobolev space with norm

\be
\|u\|_{2-\frac{2}{q},q} = \|u\|_{1,q} + \left\{\int_{\Omega}\limits\int_{\Omega}\limits \frac{|u'(x)-u'(y)|^q}{|x-y|^{1+\sigma q}}~dxdy\right\}^{\frac{1}{q}},\label{frac_norm}
\ee
where $\Omega=(0,1)$, $\sigma=1-\frac{2}{q}$.

Then it is elementary to prove
\\~\\
\textbf{Lemma A.2}~~~~$H^2((0,1))\hookrightarrow W^{2-\frac{2}{q},q}((0,1))$, {\em for $1<q<4$}.
\\~\\

Finally, let $W^{2,1}_q(Q_T)$ be the $L^q$-type parabolic Sobolev space with norm
\be
\|u\|^{(2)}_{q,Q_T}=\|u\|_q(Q_T)+\|u_x\|_q(Q_T)+\|u_{xx}\|_q(Q_T)+\|u_t\|_q(Q_T).
\ee
Then, by the fractional-order case of \cite{lady}, Ch. IV, Thm. 9.1, we have
\\~\\
\textbf{Lemma A.3}~~~~If $u$ is the solution of the mixed Dirichlet/Neumann problem for (\ref{linear_para}), and if $\|u_x\|_{\infty}$ is {\em a priori} bounded, then we have
\be
\|u\|^{(2)}_{q,Q_T}\leq c\left(\|u(\cdot,0)\|_{2-\frac{2}{q},q}(\Omega) + |u_1|\right).\label{W(2)_q_est}
\ee

Finally, by \cite{lady}, Ch. IV, top of p.343, we have
\\~\\
\textbf{Lemma A.4}~~~~For $u\in W^{2,1}_q(Q_T)$, and $q>3$, there holds the imbedding inequality
\be
|u|^{(2-\frac{3}{q})}_{Q_T}\leq c \|u\|^{(2)}_{q,Q_T}.
\ee

\section*{Appendix B. A useful inequality}

\textbf{Lemma B.1}~~~~{\em For a function $u\in H^{2+l,\frac{1}{2}(2+l)}(\Omega_T)$, $0<l<1$, such that $\Omega_T=\Omega\times[0,T]$,  $\Omega\subset\mathbb{R}$, we have}
\be
|u|^{(l+1)}_{\Omega_T}\leq C(\Omega,l)\left(T^\delta|u|^{(l+2)}_{\Omega_T} + \|u(\cdot,0)\|_{C^2(\Omega)}\right),
\ee
{\em where $\delta=\min\{l/2,(1-l)/2\}$}.

{\em Proof.}~~~~The parabolic H\"{o}lder norms in question are, when written out in full,

\bea
|u|^{(l+1)}_{\Omega_T} & := & \|u\|_{C(\Omega_T)} + \|u_x\|_{C(\Omega_T)} + \langle u_x\rangle^l_{x,\Omega_T}\nonumber\\
& + & \langle u\rangle^{(l+1)/2}_{t,\Omega_T} + \langle u_x\rangle^{l/2}_{t,\Omega_T},\nonumber\\
~ & ~ & ~ \nonumber\\
|u|^{(l+2)}_{\Omega_T} & := & \|u\|_{C(\Omega_T)} + \|u_x\|_{C(\Omega_T)} + \|u_{xx}\|_{C(\Omega_T)} + \|u_t\|_{C(\Omega_T)}\nonumber\\
& +  & \langle u_x\rangle^l_{x,\Omega_T} + \langle u_{xx}\rangle^l_{x,\Omega_T} + 
\langle  u_t\rangle^l_{x,\Omega_T}\nonumber\\
& +  & \langle u_t\rangle^{l/2}_{t,\Omega_T} + \langle u_x\rangle^{(l+1)/2}_{t,\Omega_T} +
\langle u_{xx}\rangle^{l/2}_{t,\Omega_T}.
\eea

We proceed to estimate each of the terms appearing in $|u|^{(l+1)}_{\Omega_T}$. 

First,

\be
\|u\|_{C(\Omega_T)}\leq T\|u_t\|_{C(\Omega_T)} + \|u(\cdot,0)\|_{C(\Omega)}.
\ee

Second,

\bea
\|u_x\|_{C(\Omega_T)} & \leq &\sup_{x,t}|u_x(x,t)-u_x(x,0)| + \sup_{x}|u_x(x,0)|\nonumber\\
 & \leq & T^{(l+1)/2}\langle u_x\rangle^{(l+1)/2}_{t,\Omega_T} + \|u_x(\cdot,0)\|_{C(\Omega_T)}.
\eea

Third,

\bea
\langle u_x\rangle^l_{x,\Omega_T} & = & \sup_{x,x',t}\frac{|u_x(x,t)-u_x(x',t)|}{|x-x'|^l}\nonumber\\
 & = & \sup_{x,x',t}\frac{|\int_{x'}^x u_{xx}(x'',t)~dx''|}{|x-x'|^l}\nonumber\\
 & \leq & \sup_{x,x'}|x-x'|^{1-l}\|u_{xx}\|_{C(\Omega_T)}\nonumber\\
 & = & C(\Omega,l)\|u_{xx}\|_{C(\Omega_T)}\nonumber\\
 & \leq & C\sup_{x,t}(|u_{xx}(x,t)-u_{xx}(x,0)| + |u_{xx}(x,0)|)\nonumber\\
 & \leq & CT^{l/2}\langle u_{xx}\rangle^{l/2}_{t,\Omega_T} + \|u_{xx}(\cdot,0)\|_{C(\Omega)}.
\eea

Fourth,

\bea
\langle u\rangle^{(l+1)/2}_{t,\Omega_T} & = & \sup_{x,t,t'}\frac{|u(x,t)-u(x,t')|}{|t-t'|^{(l+1)/2}}\nonumber\\
 & \leq & \|u_t\|_{C(\Omega_T)}T^{(1-l)/2}.
\eea

Fifth,

\bea
\langle u_x\rangle^{l/2}_{t,\Omega_T} & = & \sup_{x,t,t'}\frac{|u_x(x,t)-u_x(x,t')|}{|t-t'|^{l/2}}\nonumber\\
 & \leq & T^{\frac{1}{2}}\langle u_x\rangle^{(l+1)/2}_{t,\Omega}.
\eea

Putting these estimates together gives the desired result $\square$

\section*{Appendix C. Weak formulation of the multi-phase problem}

Subject to a restriction on the initial data, a weak formulation of $SP_n$ can be written down in such a way that the Rankine-Hugoniot condition becomes `hidden', thus aiding mathematical analysis. Essentially the same kind of formulation was written down many years ago by Lady\v{z}enskaya {\em et al.} \cite{lady}, and, for example, their uniqueness proof goes through without change.

For the construction that follows, we are forced to assume that 

\bd
\textrm{(C1)}\qquad\rho\leq\rho_1\quad \textrm{in low-density phases, and}\quad\rho\geq\rho_2\quad\textrm{in high-density ones}. 
\ed

Recall that on each phase we have

\be
\frac{\p\rho}{\p t} = \Delta K(\rho),
\ee
and that the Rankine-Hugoniot condition (\ref{RH}) is satisfied at each jump. Since $K(\rho_1)=K(\rho_2)$, we can define a somewhat flattened $\widetilde{K}(\rho)$ by

\be
\widetilde{K}(\rho) = \left\{
\begin{array}{ccc}
K(\rho) & : & \rho\notin[\rho_1,\rho_2]\\
K(\rho_1) (=K(\rho_2)) & : &  \rho\in[\rho_1,\rho_2],
\end{array}
\right.
\ee
such that $\widetilde{K}$ has a piecewise-smooth, monotonically-increasing inverse, which we denote by $b$; the functions $\widetilde{K}$ and $b$ are depicted in Figure 8.

\begin{figure}

\centering

\resizebox{5in}{3.5in}{\includegraphics{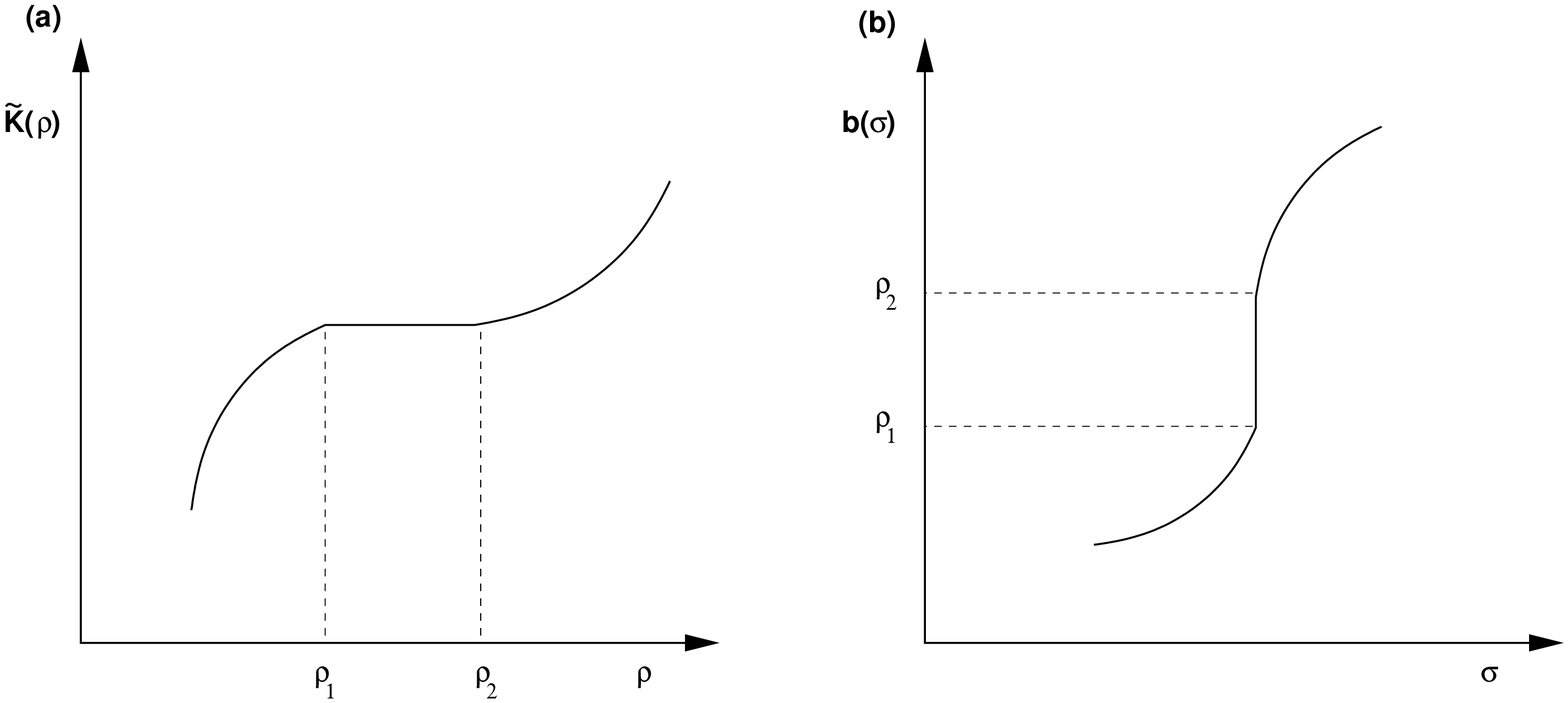}}

\caption{The functions $\widetilde{K}(\rho)$ and $b(\sigma)$.}

\end{figure}

If we introduce the new independent variable $\sigma=\widetilde{K}(\rho)$, then clearly the equation
\be
\frac{\p}{\p t}b(\sigma) = \Delta\sigma\label{cont_sigma}
\ee
is satisfied in each phase, by assumption (C1).

In terms of $\sigma$, the jump condition takes the simple form
\be
\frac{ds_i}{dt} = -\frac{(\sigma_x(s_i^+) - \sigma_x(s_i^-))}{(\rho_2-\rho_1)},\label{RH_sigma}
\ee
for $i=1,\ldots,n$.

To obtain the correct weak formulation on $Q_T = [0,1]\times[0,T]$, first note that for a smooth test function $\phi(x,t)$ such that $\phi(x,T)=0$, we have, for classical solutions of $SP_n$, and with $\rho:=b(\sigma)$,

\be
\frac{d}{dt}\int_0^1 b(\sigma)\phi~dx = \sum_{i=1}^n \dot{s}_i[\rho]_i\phi(s_i,t) + \int_0^1 \left(\phi\p_tb(\sigma) + b(\sigma)\p_t\phi\right)~dx
\ee

\be
\Longrightarrow\quad - \int_0^1 b(\sigma)\phi(x,0)~dx = \int_0^T\sum_{i=1}^n \dot{s}_i[\rho]_i\phi(s_i,t)~dt + \int_0^T\!\int_0^1\left( \phi\p_tb(\sigma) + b(\sigma)\p_t\phi\right)~dxdt,
\ee
where $[\rho]_i$ is the leap of $\rho$ at $s_i$.

For the right-hand side of (\ref{cont_sigma}), we have the weak form

\be
\int_0^T\!\int_0^1 \phi\Delta\sigma~dxdt = \int_0^T\!\int_0^1 \sigma\Delta\phi~dxdt -  \int_0^T\sum_{i=1}^n\phi(s_i,t)[\sigma_x]_i~dt - \int_0^T\sum_{i=1}^n \phi_x(s_i,t)[\sigma]_i~dt.
\ee
Thus, using (\ref{RH_sigma}), and noting that $\sigma$ is continuous at $s_i$ by construction, we arrive at
\be
\int_0^T\!\int_0^1 (b(\sigma)\p_t\phi + \sigma\Delta\phi)~dxdt + \int_0^1 b(\sigma)\phi|_{t=0}~dx = 0\label{weak}
\ee
for all smooth $\phi(x,t)$ such that $\phi(x,T)=0$, as the weak formulation of $SP_n$.

This is identical to the problem considered in Chapter V.9 of \cite{lady}, and existence and uniqueness of solutions follows by exactly the same argument used there, since one merely requires that $b(\sigma)$ be piecewise smooth and monotonically increasing. We record this result as
\\~\\
\textbf{Lemma C.1}~~~~{\em For a given bounded, continuous initial datum $\psi(x)$, equation (\ref{weak}) has a unique bounded solution $\sigma(x,t)$.}


\begin{thebibliography}{100}
\bibitem{adams} Adams, R.: Sobolev Spaces. Academic Press (1975).
\bibitem{amann} Amann, H.: Dynamic theory of quasilinear parabolic systems. III. Global existence. Math. Z. \textbf{202}, 219-250 (1989).
\bibitem{anguige} Anguige, K., Schmeiser, C.: A one-dimensional model of cell diffusion and aggregation, incorporating volume-filling and cell-to-cell adhesion. J. Math. Biol., online first.
\bibitem{lady} Lady\v{z}enskaya, O., Solonnikov, V., Ural'ceva, N.: Linear and Quasilinear Equations of Parabolic Type. AMS Translations of Mathematical Monographs, Vol.23 (1968).
\bibitem{meir} Meirmanov, A.: The Stefan Problem. de Gruyter Expositions in Mathematics (1992).
\bibitem{vii} Sun, X., Ward, M.: The Dynamics and Coarsening of Interfaces for the Viscous Cahn-Hilliard Equation in One Spatial Dimension. Stud. Appl. Math. \textbf{105}, 203-234 (2000).
\bibitem{x}  Taylor, M.: Partial Differential Equations III. Springer (1996).
\bibitem{vazquez} Vazquez, J. L.: The Porous-Medium Equation: Mathematical Theory. Oxford Science Publications (2007).
\end{thebibliography}
\end{document}